\documentclass[11pt]{amsbook}
\usepackage{amsmath}
\usepackage{amssymb}
\usepackage{amscd}
\usepackage{amsthm}

\usepackage{graphicx}
\usepackage{pstricks}
\usepackage{pst-plot}

\newtheorem{thm}{Theorem}[section]
\newtheorem{lem}[thm]{Lemma}
\newtheorem{cor}[thm]{Corollary}
\newtheorem{conj}[thm]{Conjecture}
\newtheorem{prop}[thm]{Proposition}
\newtheorem{question}[thm]{Question}
\newtheorem{hyp}[thm]{Hypothesis}
\newtheorem{claim}[thm]{Claim}

\theoremstyle{remark}

\newtheorem{rem}[thm]{Remark}

\theoremstyle{definition}
\newtheorem{defn}[thm]{Definition}

\numberwithin{equation}{section}

\allowdisplaybreaks[4]

\DeclareMathOperator{\Pic}{Pic}

\DeclareMathOperator{\Div}{Div}

\DeclareMathOperator{\supp}{supp}
\DeclareMathOperator{\Hom}{Hom}
\DeclareMathOperator{\Sym}{Sym}

\DeclareMathOperator{\Bl}{Bl}

\begin{document}

\renewcommand{\qedsymbol}{q.e.d.}
\vfuzz0.5pc
\hfuzz0.5pc 

\newcommand{\claimref}[1]{Claim \ref{#1}}
\newcommand{\thmref}[1]{Theorem \ref{#1}}
\newcommand{\propref}[1]{Proposition \ref{#1}}
\newcommand{\lemref}[1]{Lemma \ref{#1}}
\newcommand{\coref}[1]{Corollary \ref{#1}}
\newcommand{\remref}[1]{Remark \ref{#1}}
\newcommand{\conjref}[1]{Conjecture \ref{#1}}
\newcommand{\questionref}[1]{Question \ref{#1}}
\newcommand{\defnref}[1]{Definition \ref{#1}}
\newcommand{\hypref}[1]{{\it \ref{#1}\/}}
\newcommand{\Hypref}[1]{{\it Hypothesis \ref{#1}\/}}
\newcommand{\secref}[1]{Sec. \ref{#1}}
\newcommand{\ssecref}[1]{\ref{#1}}
\newcommand{\sssecref}[1]{\ref{#1}}
\newcommand{\figref}[1]{Figure \ref{#1}}

\def \red{{\mathrm{red}}}
\def \tors{{\mathrm{tors}}}
\def \EQ{\Leftrightarrow}

\def \mapright#1{\smash{\mathop{\longrightarrow}\limits^{#1}}}
\def \mapleft#1{\smash{\mathop{\longleftarrow}\limits^{#1}}}
\def \mapdown#1{\Big\downarrow\rlap{$\vcenter{\hbox{$\scriptstyle#1$}}$}}
\def \smapdown#1{\downarrow\rlap{$\vcenter{\hbox{$\scriptstyle#1$}}$}}
\def \A{{\mathbb A}}
\def \CA{{\mathcal A}}
\def \CB{{\mathcal B}}
\def \I{{\mathcal I}}
\def \J{{\mathcal J}}
\def \CO{{\mathcal O}}
\def \C{{\mathcal C}}
\def \BC{{\mathbb C}}
\def \BQ{{\mathbb Q}}
\def \m{{\mathcal M}}
\def \H{{\mathcal H}}
\def \S{{\mathcal S}}
\def \Z{{\mathcal Z}}
\def \BZ{{\mathbb Z}}
\def \W{{\mathcal W}}
\def \Y{{\mathcal Y}}
\def \T{{\mathcal T}}
\def \P{{\mathbb P}}
\def \CP{{\mathcal P}}
\def \G{{\mathbb G}}
\def \F{{\mathbb F}}
\def \BR{{\mathbb R}}
\def \D{{\mathcal D}}
\def \L{{\mathcal L}}
\def \f{{\mathcal F}}
\def \X{{\mathcal X}}
\def \E{{\mathcal E}}
\def \n{{\mathcal N}}

\def \closure#1{\overline{#1}}
\def \EQ{\Leftrightarrow}
\def \imply{\Rightarrow}
\def \isom{\cong}
\def \embed{\hookrightarrow}
\def \tensor{\mathop{\otimes}}
\def \wt#1{{\widetilde{#1}}}

\title{On Algebraic Hyperbolicity of Log Varieties}

\author{Xi Chen}

\address{Department of Mathematics\\
South Hall, Room 6607\\
University of California\\
Santa Barbara, CA 93106}
\email{xichen@math.ucsb.edu}
\date{Oct 20, 2001}
\thanks{Research supported in part by a Ky Fan Postdoctoral Fellowship}

\begin{abstract}
Kobayashi conjecture says that every holomorphic map $f: \BC\to
\P^n\backslash D$ is constant for a very general hypersurface
$D\subset \P^n$ of degree $\deg D\ge 2n+1$. As a corollary of our
main theorem, we show that $f$ is constant if $f(\BC)$ is
contained in an algebraic curve.
\end{abstract}

\maketitle

\section{Introduction and Statement of Results}\label{SEC001}

\subsection{Kobayashi Conjecture}\label{SEC001001}

The hyperbolicity in the title refers to Kobayashi hyperbolicity
\cite{K}. A complex manifold $M$ is hyperbolic in the sense of
S. Kobayashi if the hyperbolic pseudo-metric defined on $M$ is a metric.
One consequence of a complex manifold $X$ being hyperbolic is that
there does not exist nonconstant holomorphic map from $\BC$ to $X$;
actually, this is also sufficient if $X$ is compact by R. Brody
\cite{B}.

It is usually very hard to prove a complex manifold to be hyperbolic,
even for very simple ones, such as hypersurfaces and
their complements in projective spaces.
S. Kobayashi conjectured that the following is true.

\begin{conj}[Kobayashi Conjecture]\label{CONJ001}
For $n\ge 2$,
\begin{enumerate}
\item a very general hypersurface $D\subset\P^{n+1}$ of
degree $\deg D\ge 2n + 1$ is hyperbolic;
\item $\P^n \backslash D$ is hyperbolic and hyperbolically embedded in
$\P^n$ for a very general hypersurface
$D\subset \P^n$ of degree $\deg D \ge 2n + 1$.
\end{enumerate}
\end{conj}

There is an ample literature on the problem. Please see, e.g.,
\cite{D} and \cite{Z1} for surveys on the subject.
There have been some major breakthroughs
recently on the conjecture. J.P. Demailly and J. El Goul proved
the conjecture for $n = 2$ and $\deg D \ge 21$ \cite{DEG}.
Independently,
M. McQuillan also proved the conjecture for $n = 2$ and $\deg D \ge
36$ \cite{M}. Later, for the second part of the conjecture, J. El Goul improved
the bound for $\deg D$ to $15$ \cite{EG}.

The purpose of this paper is not on Kobayahsi conjecture itself but
rather on a closely related problem. Our main question is: 

\begin{question}\label{Q004}
Does $\P^n\backslash D$ contain an algebraic torus
$\BC^*\hookrightarrow \P^n\backslash D$,
embedded algebraically, for $\deg D \ge 2n+1$?
\end{question}

The answer
to this question has to be positive if Kobayashi conjecture holds
true; otherwise, there will be a nonconstant map $\BC \to \BC^*
\hookrightarrow \P^n\backslash D$ with $\BC\to\BC^*$ the covering map given by
$z\to e^z$. Indeed, this question would have to be addressed
if there were to be a proof for the full statement of the conjecture.

Simple as this problem sounds, surprisingly there is no proof for it
existing
in the literature or it may have escaped the attention of the researchers
in this field, especially algebraic geometers. In any case, except the
result of G. Xu in $n = 2$ which we will discuss later, there is no
affirmative answer to this question in general. Let us first set
ourselves a rather modest goal to prove the following.

\begin{claim}\label{CLM002}
There is no $\BC^*$ embedded algebraically in $\P^n\backslash D$ for a
very general hypersurface $D$ of sufficiently large degree.
\end{claim}

\begin{proof}
Since the statement is algebraic and generic in nature, we may
approach it via a degeneration argument. First, let us work with
some special $D$, say, $D = H_1\cup H_2 \cup ...\cup H_{2n+1}$ a
union of $2n+1$ hyperplanes in general position. To show there is no
$\BC^*$ in $\P^n\backslash D$, it is in essence to show there is no
rational curve $C\subset \P^n$ which meet $D$ (set-theoretically) at
no more than two distinct points. To achieve the least number
of intersections between a curve $C\subset \P^n$ and $D = H_1\cup H_2
\cup ...\cup H_{2n+1}$, the best one can do is to choose $C$ such that
$C$ meets $H_1, H_2, ..., H_n$ only at the point $p = H_1\cap H_2\cap
... \cap H_n$ and meets $H_{n+1}, H_{n+2}, ..., H_{2n}$ only at the point
$q = H_{n+1} \cap H_{n+2}\cap ... \cap H_{2n}$. But then $H_{2n+1}$
must meet $C$ at a third point since $p,q \not\in H_{2n+1}$ when $H_i$
are chosen to be in general position. So every curve $C\not\subset D$
meets $D$ at no less than three distinct points. More generally, if
$D$ is a union of $k$ hypersurfaces that meet properly among
themselves, i.e., any $n+1$ of them have empty intersection, then
every curve $C\not\subset D$ meets $D$ at no less than $\lceil k/n
\rceil$ distinct points.

But we are really interested in irreducible $D$'s. An obvious way to do
this is to degenerate $D$ to a union of hypersurfaces. The basic setup
is as follows.

Let $Z = \P^n\times \Delta$ and $W\subset Z$ be a pencil of
hypersurfaces of degree $d$ whose central fiber $W_0$ is a union
$D_1\cup D_2\cup ... \cup D_k$ of $k$ hypersurfaces.
Throughout the paper, we always use the notation $\Delta$ to denote the disk
parametrized by $t$.

Let $Y$ be a reduced irreducible flat family of rational
curves over $\Delta$ with the commutative diagram
\begin{equation}\label{E204}
\begin{CD}
Y @>\pi>> Z\\
@VVV @VVV\\
\Delta @>>> \Delta
\end{CD}
\end{equation}
where $\Delta\to \Delta$ is a base change and $\pi: Y\to Z$ is a
proper morphism. We assume that $\pi(Y)$ meets $W$ properly in
$Z$. Our goal is, of course, to show that $\pi(Y_t)$ meets $W_t$ at no
less than three distinct points for $t\ne 0$ if $d >> 0$. We do know
that $\pi(Y_0)$ meets $W_0$ at no less than $\lceil k/n \rceil$
distinct points. It seems that we are done as long as we take $d = k
\ge 2n + 1$ since the number of intersections $\pi(Y_t)\cap W_t$ is
obviously lower semi-continuous in $t$. However, there is a serious
flaw in our argument: even if $\pi(Y_t)$ meets $W_t$ properly,
$\pi(Y_0)$ could very well fail to meet $W_0$ properly, i.e., some
component of $\pi(Y_0)$ may be contained in $W_0 = D_1\cup D_2 \cup
... \cup D_k$.

There are more sophisticated ways we will later introduce
to deal with the situation that
$\pi(Y_0)$ and $W_0$ do not intersect properly. For the moment, let us
get around the problem using a very simple argument. 
Since $Y_t$ is rational, every component of $Y_0$ is rational.
Thus we may make the degree of $D_i$ large enough so that $D_i$ does
not contain any rational curves for $i = 1,2,...,k$. This is possible
thanks to a theorem of H. Clemens, which we will discuss in details
later. It says that a very general hypersurface $D$ of degree $\deg D \ge
2n - 1$ in $\P^n$ does not contain any rational curves. So it is
enough to take $\deg D_i \ge 2n - 1$. This, inevitably, will make
$d = \deg D$ really big. But at least we have achieved our modest goal
\ref{CLM002}.
\end{proof}

Actually, we have proved:

\begin{prop}\label{PROP009}
Let $D$ be a very general hypersurface of degree $d$ in $\P^n$. Then
every rational curve $C\not\subset D$ meets $D$ at no less than
\begin{equation}\label{E408}
\left\lceil \frac{1}{n}\left \lfloor \frac{d}{2n-1} \right \rfloor
\right \rceil
\end{equation}
distinct points. Consequently, $\P^n\backslash D$ does not contain
$\BC^*$ if $d \ge 4n^2 - 1$. 
\end{prop}

This is, however, nowhere near a satisfactory answer to our
question. The bound $4n^2 - 1$ in the above proposition is a far cry
from the conjectured $2n + 1$, especially in the light of G. Xu's
result in $n = 2$. In \cite{X2}, he proved that every irreducible
curve (not only rational ones) $C\not\subset D$ meets $D$ at no less
than $d - 2$ distinct points for a very general curve
$D$ of degree $d$ in $\P^2$. In the end of his paper, he asked whether
the following similar statement holds in higher dimensions.

\begin{question}\label{Q001}
If $D\subset \P^n$ is a very general hypersurface of degree $d$ and
$C\subset \P^n$ is an irreducible curve with $C\not\subset D$,
what is the least
number of distinct points that $C\cap D$ must contain? The expected
answer is $d - 2n + 2$.
\end{question}

An answer to the above question, of course, also
answers \questionref{Q004}.

Xu's proof in $n = 2$ is 
based upon the study of the deformation of the pair $(C, D)$ with 
$|C\cap D|$ fixed,
where $|C\cap D|$ is the number of distinct points in
the (set-theoretical) intersection $C\cap D$. Xu's method seems very
hard to generalize to higher dimensions or to underlying spaces other than
$\P^n$. An alternative approach based
upon degeneration was adopted in \cite{C1}. It can be summarized as
follows.

With $Y,Z$ and $W$ defined as before, we try to give a lower bound for
$|\pi(Y_t)\cap W_t|$ based upon the information on $|\pi(Y_0)\cap
W_0|$. Here $Z = \P^2\times\Delta$ and $W\subset Z$ is a pencil of
degree $d$ curves with reducible central fiber. As we pointed out
before, if $\pi(Y_0)$ fails to meet $W_0$ properly, we cannot say much
about $|\pi(Y_t)\cap W_t|$ in general. However, we can work around
this problem due to the fact that $\pi_* Y$ is a divisor in $Z$ here. The
intersection $\pi_* Y_t\cap W_t$ can be regarded as a member in the
linear series $\P H^0(W_t, \CO_{W_t}(\pi_* Y_t))$. There is a
well-defined notion of taking limit of a linear series, i.e.,
$\lim_{t\to 0} \P H^0(W_t, \CO_{W_t}(\pi_* Y_t))$. Correspondingly, we
may take the limit $\lim_{t\to 0} (\pi_* Y_t\cap W_t)$ as a member of
the limit linear series, which can be described explicitly in terms of
the information we have on $\pi_* Y_0$, $W_0$ and the base locus of the
pencil $W$.

This approach seems more adaptable than Xu's.
It can be easily extended to surfaces other than $\P^2$ as in
\cite{C1}. However, this line of argument still cannot be carried over
to higher dimensions, since in order for this to work, we need to treat
$\pi_* Y_t \cap W_t$ as an element in the Chow group $A_0(W_t)$; however,
with our present knowledge of Chow groups, we have no idea how to take
the limit $\lim_{t\to 0} A_0(W_t)$, even only for $W_t$ surfaces.

Note that both the above and Xu's argument did not use the
information on the genus of $Y_t$, while the fact that $Y_t$ is
rational is used in an essential way in \ref{CLM002}. In order to
answer \questionref{Q004},
it is suggested and carried out for surfaces in \cite{C2} to
bound the geometric genus $g(C)$ and $|C\cap D|$ at the same
time. More precisely, we are trying to bound the quantity
\begin{equation}\label{E409}
2 g(C) - 2 + |C\cap D|.
\end{equation}
Of course, we are done if we can show that \eqref{E409} is positive for
every reduced irreducible curve $C\not\subset D$ if $\deg D \ge 2n+1$.

Another subtle point is that instead of using the number of
the set-theoretical intersections between $C$ and $D$, we should use a
more natural notion of intersection as defined in \cite{C2}.

\begin{defn}\label{DEF007}
Let $D$ be an effective divisor on $X$ and $C\subset X$ be a reduced
irreducible curve such that $C\not\subset D$. Let $\nu:
C^{\nu}\to C\subset X$ be the normalization of $C$ and then $i_X(C,
D)$ is the number of distinct points in the set
$\nu^{-1}(D)$.
\end{defn}

Roughly, $i_X(C, D)$ is the intersection between the normalization of
$C$ and $D$. This is a more natural notion than $|C\cap D|$ since $i_X(C,
D)$ depends only on the complement $X\backslash D$ while $|C\cap D|$
also depends on the choice of the compactification of $X\backslash
D$. For example, suppose
that $C$ meets $D$ at an ordinary double point, i.e., a node $p$
of $C$. Let $X'$ be the blowup of $X$ at $p$, $D'$ be the total
transform of $D$ and $C'$ be the proper transform of $C$. Then we have
$X\backslash D\isom X'\backslash D'$ and $i_X(C, D) = i_{X'}(C', D')$ but 
$|C'\cap D'| = |C\cap D| + 1$.

So we are actually trying to bound the quantity
\begin{equation}\label{E410}
2 g(C) - 2 + i_X(C, D).
\end{equation}
Our main result is the following.

\begin{thm}\label{THM001}
Let $D\subset \P^n$ be a very general hypersurface of degree $d$. Then
\begin{equation}\label{E304}
2 g(C) - 2 + i_{\P^n}(C, D) \ge (d - 2n) \deg C
\end{equation}
for all reduced irreducible curves $C$ with $C\not\subset D$.
\end{thm}

Although we does not really answer \questionref{Q001} in the above
theorem, it does give us what we need as far as \questionref{Q004} is
concerned. Namely, we have the following corollary of \thmref{THM001}.

\begin{cor}\label{COR001}
For a very general hypersurface $D\subset\P^{n}$ of
degree $\deg D\ge 2n + 1$, $\P^{n}\backslash D$ does not contain any
algebraic torus $\BC^*$. Therefore, a holomorphic map $f: \BC\to
\P^{n}\backslash D$ is constant if $f(\BC)$ is contained in an
algebraic curve.
\end{cor}

In addition, \eqref{E304} also shows that $i_{\P^n}(C, D)$ goes up as
$\deg C$ goes up, which is not predicted by \questionref{Q001}.
Also it is sharp in the sense that for each $d\ge 2n$,
there exists a line meeting $D$ at exactly $d - 2n + 2$ distinct points
(see e.g. \cite{Z2}) and hence the equality in \eqref{E304} holds.

\subsection{Basic strategy}\label{SEC001002}

In \cite{C2}, we proved \thmref{THM001} for $n = 2$. In order to prove
it for $n > 2$, we need to generalize the techniques in \cite{C2} to
higher dimensions. This is mainly what this paper is about.

The proof of \thmref{THM001} for $n > 2$ is actually not very much
different to that for $n = 2$ or the proof for \ref{CLM002},
at least in principle.
It basically consists of two parts:
\begin{enumerate}
\item we first prove that \eqref{E304} holds when $D$ is a union of $d$
hyperplanes in general position;
\item we then prove \eqref{E304} for $D$ irreducible by
degenerating it to a union of hyperplanes.
\end{enumerate}
However, both parts of the proof are technically harder than those for
$n = 2$, especially the degeneration part.

Let us first recall some definitions and notations employed in \cite{C2}.

We have already defined the number $i_X(C, D)$ for a log pair $(X,
D)$ and a curve $C\subset X$. For technical
reasons, we will define $i_X(C, D, P)$ involving an extra term $P$
as follows. 

\begin{defn}\label{DEF002}
Let $X$ be a scheme, $D\subset X$ be a closed subscheme pure of
codimension one, $P\subset X$ be a closed subscheme of $X$ and
$C$ be a reduced irreducible curve on $X$. Let $f: \wt{X}\to X$ be the
blowup of $X$ along $P$ and let $\wt{C}$ and $\wt{D}$ be the proper transforms 
of $C$ and $D$ under $f$, respectively. If $\wt{C}$ and $\wt{D}$ exist
and meet properly, then we define
\begin{equation}\label{E101}
i_X(C, D, P) = i_\wt{X}(\wt{C}, \wt{D}). 
\end{equation}
It is also convenient to extend the definition of
$i_X(C, D, P)$ to the following situations.
\begin{enumerate}
\item Suppose that $C\subset D = \cup_{i\in I} D_i$ with $D_i$
irreducible. We define
\begin{equation}\label{E102}
i_X(C, D) = i_X(C, \cup_{i\not\in J} D_i),
\end{equation}
where $J = \{j\in I: C\subset D_j\}$.
\item Suppose that $\wt{C}$ does not exist, i.e., $C$ is contained in
the exceptional locus of $f$. Then $i_X(C, D, P)$ is defined by
\begin{equation}\label{E103}
i_X(C, D, P) = \inf_{\genfrac{}{}{0pt}{}{\Gamma\subset f^{-1}
C}{f_* \Gamma = C}} i_\wt{X}(\Gamma, \wt{D})
\end{equation}
where we consider all curves $\Gamma\subset \wt{X}$ that map
birationally onto $C$ by $f$.
\item If $C$ is nonreduced or reducible, we define
\begin{equation}\label{E104}
i_X(C, D, P) =
\sum_{\Gamma\subset C} \mu_\Gamma i_X(\Gamma, D, P),
\end{equation}
where we sum over all irreducible components $\Gamma\subset C$ and
$\mu_\Gamma$ is the multiplicity of $\Gamma$ in $C$.
\end{enumerate}
\end{defn}

A central theme of \cite{C2} is the algebraic hyperbolicity of a log
variety $(X, D)$.

\begin{defn}\label{DEF003}
Let $X\subset \P^N$ be a projective variety, $D$ be an effective
divisor on $X$ and $P\subset D$ be a closed subscheme of $X$ of
codimension at least $2$. We call $(X, D, P)$
{\it algebraically hyperbolic\/} if there exists a positive number
$\epsilon$ such that
\begin{equation}\label{E305}
2 g(C) - 2 + i_X(C, D, P) \ge \epsilon \deg C
\end{equation}
for all reduced irreducible curves $C\subset X$ with $C\not\subset
D$. And we call $(X, D)$ algebraically hyperbolic if $(X, D,
\emptyset)$ is.
\end{defn}

So \thmref{THM001} implies that $(\P^n, D)$ is algebraically
hyperbolic if $\deg D\ge 2n+1$.

Our definition of algebraic hyperbolicity for log varieties
is a natural extension of J.P. Demailly's definition for projective
varieties \cite[Chap. 2]{D}, i.e., a projective variety $X\subset\P^N$ is
called algebraically hyperbolic if $(X, \emptyset)$ is, i.e., there
exists a positive number $\epsilon$ such that
\begin{equation}\label{E306}
2g(C) - 2\ge \epsilon \deg C
\end{equation}
for all reduced irreducible curves $C\subset X$. Demailly proved
that there are no nonconstant maps from an abelian variety to an
algebraically hyperbolic variety \cite[Theorem 2.1, p. 293]{D}. Using
his argument, one can show a similar statement for log varieties.

\begin{prop}\label{PROP005}
If $(X, D)$ is algebraically hyperbolic, then there are no
nonconstant maps from a semiabelian variety to $X\backslash D$.
\end{prop}
 
In order to state our next theorem, we need to introduce the following
terms.

\begin{defn}\label{DEF004}
We call a nonempty finite
set $\{ D_i \}_{i \in I}$ of base point free (BPF) divisors on $X$
{\it an effective adjunction sequence\/}
if for any two disjoint subsets $I_1$ and $I_2$ of $I$ satisfying 
$|I_1|+|I_2| = |I| - 1$, $\sum_{i\in I_1} D_i$
is ample when restricted to $\cap_{i\in I_2} D_i$.
Here by a BPF divisor, we refers to a
general member of a BPF linear system. Note that if $I_1 = \emptyset$,
$\sum_{i\in I_1} D_i = 0$; for it to be ample on $\cap_{i\in I_2}
D_i$, it is necessary that $\dim \cap_{i\in I_2} D_i\le 0$.
\end{defn}

\begin{rem}\label{REM004}
Let $n = \dim X$. Obviously, $n+1$ BPF ample divisors on $X$ form an effective
adjunction sequence. So $n+1$ hyperplanes in $\P^n$ form an
effective adjunction sequence. If $\{ D_i\}_{i\in I}$ is an
effective adjunction sequence on $M$ and $\{ E_j \}_{j\in J}$
is an effective adjunction sequence on $N$, then
$\{ \pi_M^* D_i\}_{i\in I}\cup \{\pi_N^* E_j \}_{j\in J}$
is an effective adjunction sequence on $M\times N$,
where $\pi_M: M\times N\to M$ and $\pi_N: M\times N\to N$ are
the projections from $M\times N$ to $M$ and $N$, respectively.

We will justify the name ``effective adjunction sequence'' in
\ssecref{SEC002001}.
\end{rem}

We use $C_1(X)$ to denote the free abelian group generated by
the curves on a variety $X$ and let $N_1(X) = C_1(X)/\sim_{num}$,
where $\sim_{num}$ is the numerical equivalence.
We call $\varphi: N_1(X)\to \BR$ an additive function on $N_1(X)$ if
$\varphi\in \Hom(N_1(X), \BR)$. Of course, if $X$ is nonsingular,
$\Hom(N_1(X), \BR) \isom N^1(X)\tensor \BR$ where $N^1(X) =
\Div(X)/\sim_{num}$, i.e., every additive function $\varphi$ on
$N_1(X)$ is given by a $\BR$-divisor $D$ such that $\varphi(C)
= D\cdot C$ for all curves $C\subset X$. Let $\overline{NE}(X)\subset
N_1(X)$ be the cone of effective curves on $X$. We call $\varphi:
N_1(X)\to \BR$ a {\it subadditive\/} function on $N_1(X)$ if
\begin{equation}\label{E310}
\varphi(C_1) + \varphi(C_2) \ge \varphi(C_1 + C_2)
\end{equation}
for any $C_1, C_2\in \overline{NE}(X)$.

\begin{defn}\label{DEF001}
Let $\{ D_i \}_{i\in I}$ be a finite set of BPF divisors on
a variety $X$, $\varphi: C_1(X)\to \BR$ be
a function on $C_1(X)$ and $j_1, j_2 \in \BR$.
We say that $\{ D_i \}_{i\in I}$
satisfies the condition $\C_X(\varphi, j_1, j_2)$ on $X$ if the
following holds: for any two disjoint subsets $J_1$ and $J_2$ of $I$
satisfying $|J_1| \le j_1 $, $|J_2| \le j_2$ and 
\begin{equation}\label{E311}
\prod_{j\in J_1} D_j \ne 0,
\end{equation}
there exists a partition $I = I_1\sqcup I_2\sqcup J_2$ of
$I$ with the properties that $J_1\subset I_1$,
$\{ D_i \}_{i\in I_1}$ is an effective adjunction sequence on $X$ and
\begin{equation}\label{E312}
\sum_{i\in I_2} D_i \cdot C \ge \varphi(C)
\end{equation}
holds for all curves $C\subset X$.
\end{defn}

\begin{thm}\label{THM002}
Let $X$ be a smooth projective variety of dimension $n$,
$\{ \P \L_i \}_{i\in I}$ be a finite set of BPF linear systems on
$X$ and $D_i$ be a very general member of
$\P\L_i$ for each $i\in I$.

Let $\varphi: C_1(X) \to \BR$ be a function on $C_1(X)$ such that
$\{ D_i \}_{i\in I}$ satisfies the condition
$\C_X(\varphi, 0, n - 1)$ on $X$.

Let $F$ be a fixed effective divisor of $X$ 
and $P\subset F\cap D$ be a closed subscheme of $X$ pure of codimension
two or empty, where $D = \sum_{i\in I} D_i$. Then
\begin{equation}\label{E303}
2 g(C) - 2 + i_X(C, D, P) \ge \varphi(C)
\end{equation}
for all reduced irreducible curves $C\subset X$ and $C\not\subset F$.
\end{thm}

\begin{rem}\label{REM003}
Note that $F$ is a fixed divisor while $D_i$ varies in $\P
\L_i$. The conclusion of the theorem should be more precisely
phrased as ``for a very general choice of $D_i\in \P\L_i$ for each
$i\in I$, \eqref{E303} holds for ...''. Also note that as $D_i$
varies, $P$ varies as well. But since $P$ is a closed subscheme of
$X$ pure of codimension two and $P\subset D\cap F$, $P$ is determined
by $D_i$ up to finitely many choices. That is, if we let $F =
\sum_{j\in J} d_j F_j$ where $F_j$ are the irreducible components of
$F$, then $P$ is given by fixing a set $\CP\subset I\times J$ and
$\hat d_j \le d_j$ and setting
\begin{equation}\label{E309}
P = \bigcup_{(i,j)\in \CP} (D_i\cap \hat d_j F_j).
\end{equation}
\end{rem}

Let $X = \P^n$ and $D = \sum D_i$ be a union of $d$
hyperplanes. Then we may apply the above theorem by taking $F =
\emptyset$, and $\varphi(C) = (d-2n) \deg C$. It is easy to check
that $\{D_i\}$ satisfies the condition $\C_X(\varphi, 0, n-1)$
if $d\ge 2n$ and hence
\eqref{E304} holds for a union $D$ of $d$ hyperplanes in general
position. A natural way to go from \thmref{THM002} to \ref{THM001} is to 
degenerate an irreducible $D$ to a union of hyperplanes. 
That is what we are going to do next, albeit in a general setting.

\begin{thm}\label{THM003}
Let $X$ be a smooth projective variety,
$\{ \P \L_i \}_{i\in I}$ be a finite set of BPF linear systems on
$X$ and $D_i$ be a very general member of $\P\L_i$ for
every $i\in I$.

Let
\begin{equation}\label{E108}
I = \left(\bigsqcup_{\alpha\in\CA} I_\alpha\right) \sqcup
\left(\bigsqcup_{\beta\in \CB} I_\beta\right),
\end{equation}
be a partition of $I$, $D_{I_\gamma}$ be a very general member of
the linear system $\P \L_{I_\gamma}$ for $\gamma\in \CA\cup\CB$,
$S = \cap_{\alpha\in \CA} D_{I_\alpha}$ and $D = \sum_{\beta\in \CB}
D_{I_\beta}$, where
\begin{equation}\label{E307}
\L_J = \bigotimes_{j\in J} \L_j
\end{equation}
for $J\subset I$.

Let $\varphi: N_1(X)\to \BR$ be a subadditive function on $N_1(X)$
such that $\{ D_i \}_{i\in I}$ satisfies the condition 
$\C_X(\varphi, \dim X - 2, \dim S - 1)$. Then
\begin{equation}\label{E308}
2 g(C) - 2 + i_X(C, D) \ge \varphi(C)
\end{equation}
for all reduced curves $C\subset S$.
\end{thm}

\begin{rem}\label{REM002}
For example, we let $X = \P^n$, $\varphi(C) = (d-2n) \deg C$,
$\CA = \emptyset$, $\{D_i\}$ consist of $d$ hyperplanes
and $D\in \left|\sum D_i\right|$ be a very general hypersurface of
degree $d$. It is not hard to check that $\{ D_i \}$
satisfies the condition $\C_X(\varphi, n-2, n-1)$ if $d\ge 2n$.
Therefore, \eqref{E304} holds for $d\ge 2n$.
It is easy to check that \eqref{E304} holds for $d < 2n$ as well.
Therefore, \thmref{THM001} follows.
\end{rem}

We do not require $D$ to be irreducible in the above
theorem. Therefore, \thmref{THM001} holds for $D$ reducible as well. 

\begin{cor}\label{COR002}
Let $D = \cup D_k\subset \P^n$, where $D_k\subset\P^n$ is a very
general hypersurface of degree $d_k$. Then
\begin{equation}\label{E314}
2 g(C) - 2 + i_{\P^n}(C, D) \ge (d - 2n) \deg C
\end{equation}
for all reduced curves $C\subset \P^n$ 
where $d = \sum d_k$.
\end{cor}

The above corollary is useful since the hyperbolicity of
$\P^n\backslash D$ has been studied for reducible $D$'s as well
(see e.g. \cite{G1}, \cite{G2} and \cite{DSW}).

As another corollary of \thmref{THM003}, we
take $X = \P^{n_1}\times \P^{n_2}\times ... \times \P^{n_k}$. Let
$H_j$ be the pullback of the hyperplane divisors under
the projections $X\to \P^{n_j}$, for $j = 1,2,...,k$.
We let
$\CA = \emptyset$, $\{D_i\}$ consist of $d_1$
divisors in $|H_1|$, $d_2$ divisors in $|H_2|$, ... and $d_k$ divisors
in $|H_k|$, $D$ be a very general member of $|\sum D_i|$
and $\varphi(C) = \epsilon (H_1 + H_2 + ... + H_k) C$, where
\begin{equation}\label{E200}
\epsilon = \min_{1\le j\le k} (d_j - n_j - n)
\end{equation}
and $n = \sum n_j$. It is not hard to check that $\{D_i\}$ satisfies
the condition $\C_X(\varphi, n-2, n-1)$ when $\epsilon \ge 0$.
So we have the
following corollary.

\begin{cor}\label{COR003}
Let $X = \P^{n_1}\times \P^{n_2}\times ... \times \P^{n_k}$ and
$D\subset X$ be a very general hypersurface of type $(d_1, d_2,
..., d_k)$. Then
\begin{equation}\label{E313}
2 g(C) - 2 + i_X(C, D) \ge  \epsilon \deg C
\end{equation}
for all reduced curves $C\subset X$, where $\epsilon$ is given by
\eqref{E200} and $\deg C = (H_1 + H_2 + ... + H_k) C$.
\end{cor}

So $(\P^{n_1}\times \P^{n_2}\times ... \times \P^{n_k}, D)$ is
algebraically hyperbolic if $d_j > n_j + n$ for $j = 1,2,...,k$ and
it is reasonable to expect $(\P^{n_1}\times \P^{n_2}\times ... \times
\P^{n_k})\backslash D$ to be hyperbolic for such $D$, which can be
regarded as a generalized Kobayashi conjecture. This result is again
sharp in the sense that the equality can be achieved in
\eqref{E313} when $\epsilon\ge 0$.

When $\CA \ne\emptyset$, $C$ is contained in the complete intersection
$S$. Although the introduction of the case $C\subset S$
is mainly for the purpose of induction, it is quite
interesting in its own right. For example, we let $X = \P^n$,
$\CB = \emptyset$, $\{D_i\}$ consist of $d$ hyperplanes and
$D_{I_\alpha}$ be a very general hypersurface of degree $d_\alpha$ for
$\alpha\in \CA$,
where $d = \sum d_\alpha$. It is easy to check that $\{D_i\}$
satisfies $\C_X(\varphi, n - 2, n - a - 1)$ with
$\varphi(C) = (d + a - 2n) \deg C$ for $d + a \ge 2n$, where
$a = |\CA|$. So we obtain the following corollary.

\begin{cor}\label{COR004}
Let $D_k\subset \P^n$ be a very general hypersurface of degree
$d_k$ for $k=1,2,...,a$.
\begin{equation}\label{E315}
2 g(C) - 2 \ge \left(\sum_{k=1}^a d_k + a - 2n\right) \deg C
\end{equation}
for all reduced curves $C \subset S = \cap D_k $.
\end{cor}

This is a well-known result due to Clemens \cite{Cl} for hypersurfaces and
to Ein \cite{E1} for complete intersections. See also \cite{E2}, \cite{V},
\cite{CLR}, \cite{C-L} and \cite{X1}.

As before, we may formulate a generalization of \coref{COR004} in
$X = \P^{n_1} \times \P^{n_2}\times ...\times \P^{n_k}$.

\begin{cor}\label{COR005}
Let $D_i\subset \P^{n_1} \times \P^{n_2}\times ...\times \P^{n_k}$ be
a very general hypersurface of type $(d_{i1}, d_{i2}, ..., d_{ik})$
for $i=1,2,...,a$. Then
\begin{equation}\label{E316}
2g(C) - 2\ge \min_{1\le j\le k} \left(\sum_{i=1}^a d_{ij} + a - n_j -
n\right)\deg C
\end{equation}
for all reduced curves $C\subset \cap D_i$,
where $n = \sum n_j$ and $\deg C = (H_1 + H_2 + ... + H_k) C$.
\end{cor}

Similar corollaries hold in other homogeneous spaces, e.g.,
Grassmanians, which we will not formulate here.

\subsection{Some questions}\label{SEC001003}

It is interesting to notice that although \coref{COR002} and \ref{COR003}
are sharp, \coref{COR004} and \ref{COR005} are most likely not. The latter
two imply that the corresponding complete intersections are
algebraically hyperbolic if the RHS's of \eqref{E315} and
\eqref{E316} are positive. However, there are no obvious
reasons why they should not be algebraically hyperbolic if the RHS's
are zero. Actually, we expect they are and some cases have already been
proved in \cite{C2}.

\begin{question}\label{Q002}
Under the hypothesis of \coref{COR005}, is \(\cap D_i\) algebraically
hyperbolic if
\begin{equation}\label{E300}
\min_{1\le j\le k} \left(\sum_{i=1}^m d_{ij} + m - n_j -
n\right) = 0?
\end{equation}
Especially, it is interesting to know whether a very general
hypersurface of degree $d = 2n-1$ in $\P^n$ is algebraically
hyperbolic.
\end{question}

There are discussions on the significance of this question 
in \cite{C2}.

Alternatively, one may ask what is the smallest $d$ such that a very
general hypersurface $D$ of degree $d$ in $\P^n$ is algebraically
hyperbolic. If $d = 2n-3$, $D$ contains lines (see
e.g. \cite{H}) and hence cannot be algebraically hyperbolic.
So the only degrees in doubt are $d = 2n-2$ and $2n-1$.
We believe that $D$ is algebraically hyperbolic if $d = 2n - 1$ while it
is not if $d = 2n - 2$. While we have some evidences for our
conjecture on $d = 2n-1$, our assertion on $d = 2n - 2$ is pure
speculation. But even for $d = 2n-1$, we do not know the answer
even for $n = 3$. That is, we do not even know whether a very general
quintic surface in $\P^3$ is algebraically hyperbolic or not, although
some similar statements for complete intersections have been proved in
\cite{C2}.

Despite the fact we have proved that some classes of quasi-projective
varieties, which are conjectured to be hyperbolic, are algebraically
hyperbolic, we still do not know whether ``hyperbolic'' implies
``algebraically hyperbolic''.

\begin{question}\label{Q003}
If $X\backslash D$ is hyperbolic (and hyperbolically embedded), does
this imply that $(X, D)$ is algebraically hyperbolic?
\end{question}

This is known for $D = \emptyset$ by Demailly \cite[Theorem 2.1,
p. 293]{D}. However Demailly's argument for the projective case does
not carry over to the quasi-projective case since the integral of the
curvature of the hyperbolic metric on an affine curve is not finite.
However, we do have the following.

\begin{prop}\label{PROP006}
Suppose that $(X, D, T_X)$ has a logarithmic
$k$-jet metric with strictly negative curvature for some $k$. Then
$(X, D)$ is algebraically hyperbolic.
\end{prop}

Please see \cite{D} for the definitions of jet differentials and jet
metrics and \cite{D-L} for those of logarithmic jet differentials and
logarithmic jet metrics.

The proof of the above statement follows the same argument of Theorem 8.1
in \cite{D}, where we only have to change from ``jets'' to
``logarithmic jets'' and from $T_{\bar C}$ to
$T_{\bar C}(-\log \nu^{-1} D)$. It is conjectured by Demailly that
$X$ is hyperbolic if and only if $(X, T_X)$ has strictly negative
$k$-jet curvature for $k$ large enough \cite[Conjecture 7.13,
p. 324]{D}. We may likewise form a conjecture for log varieties:
$X\backslash D$ is hyperbolic if and only if $(X, D, T_X)$ has
strictly negative logarithmic $k$-jet curvature.

\subsection{Conventions}

\begin{enumerate}
\item Throughout the paper, we will work exclusively over $\BC$.
\item By a variety $X$ being very general, we mean that $X$ lies on
the corresponding parameter space (Hilbert scheme or
moduli space) with countably many proper closed subschemes removed. So
the notion of being very general relies on the fact that the base
field $\BC$ we work with is uncountable.
\end{enumerate}

\section{Proof of \thmref{THM002}}\label{SEC002}

\subsection{Preliminaries}\label{SEC002001}

An effective adjunction sequence has the following properties.

\begin{lem}\label{LEM002}
Let $f: Y\to X$ be a morphism between two projective varieties $X$
and $Y$. Suppose that $\{ D_i \}_{i\in I}$
is an effective adjunction sequence on $X$. Let
$\wt{D}_i = f^* D_i$ for $i\in I$.
Then the following holds.
\begin{enumerate}
\item If $f$ is finite over $f(Y)$,
$\{ \wt{D}_i \}_{i\in I}$ is an effective adjunction sequence on $Y$.
\item If $f$ is finite over $f(Y)$,
for each $J\subsetneq I$, $\{ \wt{D}_i: i\in I\backslash J \}$
is an effective adjunction sequence when restricted to
$\wt{D}_J = \cap_{j\in J} \wt{D}_j$.
\item If $Y$ is smooth and $f$ is generically finite over $f(Y)$,
$K_Y + \sum_{i\in I} \wt{D}_i$ is effective.
\end{enumerate}
\end{lem}

The last statement is the reason we use the term
``effective adjunction sequence''.

\begin{proof}[Proof of \lemref{LEM002}]
The first two statements are obvious and we will leave their proofs to
the readers.

We will prove the effectiveness of $K_Y + \sum_{i\in I} \wt{D}_i$
by repeatedly cutting $f(Y)$ by $D_i$'s and
applying Kawamata-Viehweg vanishing theorem.

Obviously,
$\wt{D}_i$ is BPF on $Y$. By Bertini's theorem, we may assume that
$\wt{D}_i$ is smooth for each $i$. Since $\{D_i\}_{i\in I}$ is an
effective adjunction sequence, $\sum_{i\ne \alpha} D_i$ is ample for
any $\alpha\in I$.
So at least one of $\wt{D}_i$ is
nonempty. Let us fix $\alpha\in I$ such that 
$D_\alpha\cap f(Y)\ne \emptyset$.

Since $\sum_{i\ne \alpha} D_i$ is ample,
$\sum_{i\ne \alpha} \wt{D}_i$ is big and nef on $Y$. 
Then
\begin{equation}\label{E320}
h^1(K_Y + \sum_{i\ne \alpha} \wt{D}_i) = 0
\end{equation}
by Kawamata-Viehweg vanishing theorem. Therefore, we have the surjection
\begin{equation}\label{E321}
H^0(K_Y + \sum_{i\in I} \wt{D}_i) \twoheadrightarrow
H^0(\CO_{\wt{D}_\alpha}(K_{\wt{D}_\alpha} + \sum_{i\ne \alpha} \wt{D}_i)).
\end{equation}
Therefore, to show that $K_Y + \sum \wt{D}_i$
is effective on $Y$, it suffices to show that
$K_{\wt{D}_\alpha} + \sum_{i\ne \alpha} \wt{D}_i$ is effective on
$\wt{D}_\alpha$. Since $\{D_i: i \ne \alpha\}$ is an effective
adjunction sequence when restricted to $D_\alpha\cap f(Y)$,
$\sum_{i\ne \alpha} D_i$ is ample when restricted
to $D_\alpha\cap f(Y)$ and
hence there is at least one $D_\beta\in \{D_i: i\ne \alpha\}$ such
that $D_\alpha\cap D_\beta\cap f(Y)\ne \emptyset$. 
Continue cutting $\wt{D}_\alpha$ by
$\wt{D}_\beta$ and we obtain
\begin{equation}\label{E322}
\begin{split}
H^0(K_Y + \sum_{i\in I} \wt{D}_i) &\twoheadrightarrow
H^0(\CO_{\wt{D}_\alpha}(K_{\wt{D}_1} + \sum_{i\ne \alpha} \wt{D}_i))\\
& \twoheadrightarrow
H^0(\CO_{\wt{D}_{\alpha\beta}}(K_{\wt{D}_{\alpha\beta}} +
\sum_{i\ne\alpha,\beta} \wt{D}_i))
\end{split}
\end{equation}
where $\wt{D}_{\alpha\beta} = \wt{D}_\alpha\cap \wt{D}_\beta$.
We will carry on this argument and eventually end up with
\begin{equation}\label{E323}
H^0(K_Y + \sum_{i\in I} \wt{D}_i) \twoheadrightarrow
H^0(\CO_{\wt{D}_J}(K_{\wt{D}_J} + \sum_{i\not\in J} \wt{D}_i))
\end{equation}
where $J\subset I$ is an index subset such that $|J| = \dim Y - 1$
and
\begin{equation}\label{E326}
f(Y)\cdot \prod_{j\in J} D_j \ne 0.
\end{equation}
So it remains to justify that
$K_{\wt{D}_J} + \sum_{i\not\in J} \wt{D}_i$ is
effective. Since $\wt{D}_J$ is a curve, it suffices to show
that $\wt{D}_J \cdot \sum_{i\not\in J} \wt{D}_i \ge 2$, i.e.,
\begin{equation}\label{E324}
f(Y) \cdot \prod_{j\in J} D_j\cdot \sum_{i\not\in J} D_i \ge 2.
\end{equation}
If \eqref{E324} fails, then at most one of $f(Y) \cdot D_i \cdot
\prod_{j\in J} D_j$ is positive and the rest are zeroes for
$i\not\in J$. Let $D_\gamma\in \{D_i: i\not\in J\}$ be the divisor
such that $f(Y)
\cdot D_i \cdot \prod_{j\in J} D_j = 0$ for each $i\not\in
J\cup\{\gamma\}$. Then
\begin{equation}\label{E325}
f(Y) \cdot \prod_{j\in J} D_j\cdot \sum_{i\not\in
J\cup\{\gamma\}} D_i = 0.
\end{equation}
But \eqref{E325} contradicts the fact that $\sum_{i\not\in
J\cup\{\gamma\}}
D_i$ is ample when restricted to $\cap_{j\in J} D_j$. This finishes
the proof of \eqref{E324} and hence the effectiveness of
$K_Y + \sum_{i\in I} \wt{D}_i$.
\end{proof}

Regarding the condition $\C_X(\varphi, j_1, j_2)$, we have the
following 
observations.

\begin{lem}\label{LEM001}
If $\{D_i\}_{i\in I}$ satisfies the condition $\C_X(\varphi, j_1, j_2)$
on $X$, then the following holds.
\begin{enumerate}
\item Any $\{D_k\}_{k\in K}\supset \{D_i\}_{i\in I}$
satisfies $\C_Y(\varphi, j_1', j_2')$ on $Y$
for any $j_1' \le j_1$, $j_2' \le j_2$ and closed subscheme $Y\subset X$.
\item For any two disjoint subsets $\n_1$ and $\n_2$
of $I$ satisfying $|\n_1| \le j_1$ and $|\n_2| \le j_2$,
$\{D_i: i\in I\backslash (\n_1\cup \n_2)\}$ satisfies
$\C_Y(\varphi, j_1 - |\n_1|, j_2 - |\n_2|)$
when restricted to $Y = \cap_{i\in \n_1} D_i$.
\item Let $Y$ be a smooth projective variety and
$f: Y\to X$ be a proper morphism which is generically finite over
its image. Suppose that $j_1\ge 0$.
Then for any $J\subset I$ satisfying $|J| \le j_2$, the
curves $C\subset Y$ satisfying
\begin{equation}\label{E331}
\left(K_Y + \sum_{i\not\in J} f^* D_i\right) C < \varphi(f_* C)
\end{equation}
cannot cover $Y$.
\end{enumerate}
\end{lem}

\begin{proof}
The first statement is obvious.

For (2), let assume that $Y\ne \emptyset$; otherwise, there is nothing
to prove.
Let $J_1$ and $J_2$ be two disjoint subsets of $I\backslash (\n_1\cup
\n_2)$ satisfying $|J_1|\le j_1 - |\n_1|$ and $|J_2| \le j_2 - |\n_2|$. Since
$\{D_i\}_{i\in I}$ satisfies $\C_X(\varphi, j_1, j_2)$, there exists a
partition $I = I_1\sqcup I_2\sqcup (J_2\sqcup \n_2)$ such that
$J_1\cup \n_1\subset I_1$, $\{D_i\}_{i\in I_1}$ is an effective
adjunction sequence on $X$ and \eqref{E312} holds.
By \lemref{LEM002}, $\{D_i: i\in I_1\backslash \n_1\}$ is an effective
adjunction sequence on $Y$. Then the partition $I\backslash (\n_1\cup
\n_2) = (I_1\backslash \n_1) \sqcup I_2 \sqcup J_2$ has all the
required properties in order
for $\{D_i: i\in I\backslash (\n_1\cup \n_2)\}$
to satisfy $\C_Y(\varphi, j_1 - |\n_1|, j_2 - |\n_2|)$ on $Y$.
 
For (3), 
since $\{D_i\}_{i\in I}$ satisfies $\C_X(\varphi, j_1, j_2)$
and $|J| \le j_2$, there exists a
partition $I = I_1\sqcup I_2 \sqcup J$ of $I$ such that 
$\{D_i\}_{i\in I_1}$ is an effective adjunction sequence on $X$ and
\begin{equation}\label{E332}
\sum_{i\in I_2} f^* D_i \cdot C\ge \varphi(f_* C).
\end{equation}
If
\begin{equation}\label{E333}
\left(K_Y + \sum_{i\in I_1} f^* D_i\right) C \ge 0
\end{equation}
then
\begin{equation}\label{E334}
\begin{split}
\left(K_Y + \sum_{i\not\in J} f^* D_i\right) C
&= \left(K_Y + \sum_{i\in I_1} f^* D_i\right) C + \sum_{i\in I_2} f^*
D_i \cdot C\\
&\ge \varphi(f_* C).
\end{split}
\end{equation}
Therefore, if $C$ satisfies \eqref{E331}, we necessarily have
\begin{equation}\label{E335}
\left(K_Y + \sum_{i\in I_1} f^* D_i\right) C < 0.
\end{equation}
So such $C$ must be contained in the base locus of the linear series
$|K_Y + \sum_{i\in I_1} f^* D_i|$. By \lemref{LEM002}, $K_Y +
\sum_{i\in I_1} f^* D_i$ is effective and hence $C$ is contained in
a fixed proper closed subscheme of $Y$. Consequently, such curves
cannot cover $Y$. 
\end{proof}

\subsection{Deformation lemmas}\label{SEC002002}

Let us first recall how a similar statement of \thmref{THM002} was proved
in \cite[Theorem 1.8]{C2}. To explain the ideas of \cite{C2} in a
nutshell, we will work out an example by showing how to prove
\eqref{E304} for $D$ a union of five lines of $\P^2$.

Let $D\subset \P^2$ be a union of five lines in general position. We
claim that
\begin{equation}\label{E317}
2 g(C) - 2 + i_{\P^2}(C, D) \ge \deg C
\end{equation}
for any reduced irreducible curves $C\subset\P^2$ with
$C\not\subset D$. This was proved in two steps in \cite{C2}.

First, we prove the following: if $D\subset\P^2$ is a union of four
lines with no three lines passing through the same point, then there
are at most countably many reduced irreducible curves $C\subset\P^2$
violating \eqref{E317}.
This was proved in \cite{C2} using the results on the deformation of
maps with tangency conditions (see e.g. \cite[Chap 3, Sec
B]{H-M}). However, we find it a little cumbersome to generalize this
deformational argument to higher dimensions. So a slightly different
approach will be taken upon here (see \lemref{LEM003} below).

Second, after we proved \eqref{E317} with only countably many
exceptions for $D$ a union of four lines, the following argument
``eliminates'' these exceptions once we have an extra line. Let $D =
\cup_{i=1}^5 L_i$ be a union of five lines in general position and
let $D' = \cup_{i=1}^4 L_i$. If $C$ satisfies \eqref{E317},
we are done. Otherwise, we necessarily have
\begin{equation}\label{E318}
2 g(C) - 2 + i_{\P^2}(C, D') < \deg C.
\end{equation}
Since there are only countably many curves $C$ satisfying
\eqref{E318}, by Bertini, $L_5$ meets $C$ transversely and
$D'\cap L_5\cap C = \emptyset$. Actually, we may choose $D'$ to be the
union of any four lines of $D$ and run the same argument as
above. Eventually, we conclude that either \eqref{E317} holds or
$C$ meets $D$ transversely. If it is the latter case,
$i_{\P^2}(C, D) = 5 \deg C$ and \eqref{E317} is trivially
satisfied.

It turns out that things become more complicated if we go up in
dimensions. Let us first introduce the following definition (see also
\cite{CLR}).

\begin{defn}\label{DEF005}
Let $X$ be a scheme. We say that the subschemes $C\subset X$ with
a certain property $P$ are at most {\it $r$-filling\/} if there exists
a union $W\subset X$ of countably many closed subschemes
of $X$ of dimension $r$ such that every $C$
with property $P$ lies in $W$.
\end{defn}

Let $D\subset \P^n$ be a union of $2n + 1$ hyperplanes in general
position. We want to generalize the above argument to show that
\begin{equation}\label{E219}
2 g(C) - 2 + i_{\P^n}(C, D) \ge \deg C
\end{equation}
for any reduced irreducible curves $C\subset\P^n$ with
$C\not\subset D$. This is done in three steps.
\begin{description}
\item[Step 1] We show that the curves $C$ violating
\eqref{E219} are at most $(n-1)$-filling in $\P^n$ if $D$
is a union of $n + 2$ hyperplanes.
\item[Step 2] We show that the curves $C$ violating \eqref{E219} are
at most $(n - r)$-filling if $D$ is a union of $n + r + 1$
hyperplanes for $2\le r \le n - 1$.
\item[Step 3] We show that there are no curves $C$ violating
\eqref{E219} if $D$ is a union of $2n + 1$ hyperplanes.
\end{description}

Of course, things become more technical when we work in the
general setting of \thmref{THM002}. But our basic approach is still the same.

The following two lemmas conclude the first step of our proof of
\thmref{THM002}. These are deformation-theoretical results in essence.

Let us first recall the following definitions.

We say a variety $X$ is of {\it normal crossing\/} (has
{\it normal crossing\/}, has a {\it normal-crossing} singularity, etc.)
at a point $p\in X$ if $X$ is locally
analytically given by $\{ (z_1, z_2, ..., z_n): z_1 z_2 ... z_m =
0\}$ at $p$ for some $m\le n$. A variety is of
{\it normal crossing\/} (has {\it normal crossing\/},
has only {\it normal-crossing\/} singularities, etc.)
if it is of normal crossing everywhere. And a variety $X$ has
{\it simple normal crossing\/} if each of its irreducible components
is smooth in addition that $X$ has normal crossing.

\begin{lem}\label{LEM003}
Let $X$ be a smooth projective variety of dimension $n$ and $D$
be a reduced effective divisor on $X$. Then 
the reduced irreducible curves $C\subset X$ with the following
properties are at most $(n-1)$-filling:
$C\not\subset D$, $C$ meets $D$ only at points where
$D$ has normal crossing and
\begin{equation}\label{E341}
2 g(C) - 2 + i_X(C, D) < (K_X + D) C.
\end{equation}
\end{lem}

\begin{proof}
Let us first prove the lemma assuming that $C$ only meets $D$ at
smooth points of $D$.

Let $M$ be the variety parameterizing such curves. It does not
really matter how we construct $M$. For example,
we may realize $M$ as a subvariety of the moduli space of stable
maps to $X$ with marked points. A point in $M$ is $(f, q_1, q_2,
..., q_l)$, where $f: C^\nu\to X$ is a stable map such that $C =
f(C^\nu)$ has the required properties and $f^* D = m_1 q_1 + m_2 q_2
+ ... + m_l q_l$ with $l = i_X(C, D)$. If we construct $M$ this
way, we have got a natural compactification of $M$ and a universal
family $\pi: Y\to M$ over $M$. Compactify $M$ and
let us assume that $Y$ and $M$ are projective.

There is a natural map $f: Y \to X$. The statement of the lemma is
equivalent to that $f$ is not dominant.
Let us assume that $f$ is dominant.

Without the loss of generality, let us assume that $M$ is
irreducible; otherwise, we replace $M$ by one of its irreducible component
$M_0$ such that $\pi^{-1}(M_0)$ dominates $X$. We may also
assume that $\dim M = \dim X - 1$; if not, we may replace $M$ by a
subvariety $M_0\subset M$ satisfying that $\dim M_0 = \dim X - 1$ and
$\pi^{-1}(M_0)$ dominates $X$. So we may assume that $f: Y\to X$
is a surjective generically finite map. Let $Q_1, Q_2, ..., Q_l$ be
the $l$ sections of $\pi: Y\to M$ which restrict to the $l$ marked
points on a general fiber of $\pi$.

After we desingularize $M$ and $Y$, let us assume that $M$
and $Y$ are nonsingular. Since $f$ is surjective, we have the
exact sequence
\begin{equation}\label{E342}
0\xrightarrow{} f^* \Omega_X \xrightarrow{} \Omega_Y \xrightarrow{}
\Omega_{Y/X} \xrightarrow{} 0.
\end{equation}
Hence $c_1(\Omega_Y) = c_1(f^* \Omega_X) + c_1(\Omega_{Y/X})$.
Since $f$ is generically finite, $\Omega_{Y/X}$ is a torsion
sheaf on $Y$. A local computation shows that $Q_k \subset
\supp(\Omega_{Y/X})$ and $\Omega_{Y/X}$ has length at least $m_k -
1$ along $Q_k$ for $k=1,2,...,l$. Therefore, the divisor
\begin{equation}\label{E343}
\begin{split}
&\quad c_1(\Omega_Y) -
c_1(f^* \Omega_X) - \sum_{k=1}^l (m_k-1) Q_k\\
& = K_Y - f^* K_X -
\sum_{k=1}^l (m_k-1) Q_k
\end{split}
\end{equation}
is effective. And its restriction to a general fiber $C^\nu$
of $\pi: Y\to M$ is effective, too. Obviously, $K_Y |_{C^\nu} =
K_{C^\nu}$ and
\(\sum_{k=1}^l (m_k-1) Q_k |_{C^\nu} = \sum_{k=1}^l (m_k-1) q_k = f^* D
- \sum_{k=1}^l q_k\).
Therefore,
\begin{equation}\label{E344}
K_{C^\nu} - f^* (K_X + D) + \sum_{k=1}^l q_k
\end{equation}
is effective, which implies $2 g(C) - 2 + i_X(C, D)\ge (K_X
+ D) C$. This is a contradiction.

Now let us handle the case that $C$ meets $X$ at singular points of $D$.
Let $[C]$ be a general point of $M$. Suppose that $C$ meets $D$
at $p$ where $D$ has normal crossing.
Let $D_{sing}$ be the singular locus of $D$ and
$P$ be the irreducible component of $D_{sing}$ containing $p$. 
We blow up $X$ along $P$ and let $f: \wt{X}\to X$ be the corresponding
map.
Let $\wt{C}$ be the proper
transform of $C$ under $f$ and $\wt{D} = \supp(f^* D)$ be the
reduced total transform of $D$. It is easy to check that 
\begin{gather}
\label{E345}
g(\wt{C}) = g(C),\\
\label{E346}
i_\wt{X}(\wt{C}, \wt{D}) = i_X(C, D),\text{ and}\\
\label{E347}
K_\wt{X} + \wt{D} = f^*(K_X + D).
\end{gather}
If $\wt{C}$ meets $\wt{D}$ only at
nonsingular points of $\wt{D}$ over $p$, we are done.
If not, $\wt{C}$ will meet $\wt{D}$ at a point $p'$ on the exceptional
divisor $E$ of $f$, where $\wt{D}$ again has normal crossing.
Then we continue to blow up
$\wt{X}$ along the irreducible component of $\wt{D}_{sing}$ containing
$p'$. This process will eventually
end and \eqref{E345}-\eqref{E347} always hold during the
process.
After a finite sequence of blowups over $p$, let $\wt{X}$ be
the resulting variety, $\wt{C}$ be the proper transform of $C$ and
$\wt{D}$ be the reduced total transform of $D$. Then $\wt{C}$
will eventually meet $\wt{D}$ only at nonsingular points of
$\wt{D}$ over $p$; meanwhile, we always have
\begin{equation}\label{E348}
\begin{split}
&\quad 2 g(\wt{C}) - 2 + i_\wt{X}(\wt{C}, \wt{D}) - (K_\wt{X} + \wt{D})
\wt{C}\\
&= 2 g(C) - 2 + i_X(C, D) - (K_X + D) C
\end{split}
\end{equation}
by \eqref{E345}-\eqref{E347}. Do this for every $p\in C\cap D$
where $D$ is singular and we reduce the lemma to the case that $C$
meets $D$ only at smooth points of $D$, which we have already proved.
\end{proof}

\begin{rem}\label{REM005}
It is essential that $C$ only meets $D$ at normal
crossing singularities of $D$;
otherwise, the above lemma may not hold if worse
singularities present themselves in $C\cap D$ due to the fact that \eqref{E347}
may not hold if we blow up along $P$ where $D$ does not have normal
crossing. For example, 
let us consider the lines $L$ meeting a cubic curve in $\P^2$ at only
one point. We would expect that there are only finitely many such
lines according to the lemma. 
However, this is simply false if the cubic has a triple point,
i.e., it is the union of three lines meeting at a point $p$ for
which $L$ could be any line passing through $p$; 
on the other hand, it is fine for a nodal cubic for which $L$ could
only be one of the two tangent lines at the node or one of the three flex
lines.
%
%
%
%
\end{rem}

\begin{lem}\label{LEM004}
Let $X$ be a smooth projective variety of dimension $n$, $\{ \P
\L_i\}_{i\in I} $ be a finite set of BPF linear systems on
$X$ and $D_i$ be a general member of $\P\L_i$ for each $i\in I$.
Let $F$ be a
fixed effective divisor of $X$ and $P\subset F\cap D$ be a closed
subscheme of $X$ pure of codimension two or empty, where
$D = \sum_{i\in I} D_i$.
Then the reduced irreducible curves $C\not\subset D$ satisfying
\begin{equation}\label{E351}
2 g(C) - 2 + i_X(C, D, P) < (K_X + D) C
\end{equation}
are at most $(n-1)$-filling in $X$.
\end{lem}

There is a similar statement for surfaces in \cite[Lemma 2.1]{C2}.
The proofs of the two lemmas are very close. Basically, we will try
reduce the lemma to the case $P = \emptyset$ and $C$ meets $D$
only along the nonsingular locus of $D$ so that we may apply
\lemref{LEM003}. However, there
are some technical issues presenting themselves in high dimensions.
In the case of surfaces, $P$ is a finite set of points and we assume
$P$ to be reduced in \cite[Lemma 2.1]{C2}; hence the blowup
$\wt{X}$ of $X$ along $P$ is smooth and it is obvious that 
\begin{equation}\label{E352}
K_{\wt{X}} + f_*^{-1} D = f^*(K_X + D),
\end{equation}
where $f_*^{-1} D$ is the proper transform of $D$ under the
blowup $f: \wt{X}\to X$. Based on this observation, we reduced the
lemma to the case $P = \emptyset$ in \cite[Lemma 2.1]{C2}. However,
things are not that easy in high dimensions: $P$ could be singular;
as a consequence, the blowup $\wt{X}$ of $X$ along $P$ could
very well be singular and we no longer have \eqref{E352}. Actually,
$K_{\wt{X}}$ will even fail to be $\BQ$-Cartier if $P$ is really
``bad''. To overcome these obstacles, we need first to ``remove''
the bad singularities of $P$. This is achieved by making $F$ into
an effective divisor with simple-normal-crossing support.

\begin{proof}[Proof of \lemref{LEM004}]
It is a well-known fact
that there exists a series of blowups $f: \wt{X}\to X$ of
$X$ with smooth centers such that $f^* F$ is supported on an
effective divisor of simple normal crossing. Let
$E_f\subset \wt{X}$ be the exceptional locus of $f$,
$\wt{F} = f^*
F$, $\wt{D} = f^* D$, $\wt{P} = f^{-1}(P)$ and $\wt{C} =
f_*^{-1} C$ be the proper transform of $C$ under $f$, where we
assume that $C\subset X$ is a reduced irreducible curve such that
$C\not\subset f(E_f)$ (those curves
contained in $f(E_f)$ are at most $(n-2)$-filling anyway).
We claim that
\begin{gather}
\label{E361}
g(C) = g(\wt{C}),\\
\label{E362}
(K_X + D) C \le (K_{\wt{X}} + \wt{D}) \wt{C}, \text{ and}\\
\label{E363} 
i_X(C, D, P) \ge i_\wt{X}(\wt{C}, \wt{D}, \wt{P}) 
\end{gather}
where \eqref{E361} is obvious and \eqref{E362} follows directly
from the fact that $\wt{X}$ is the blowup of $X$ with smooth
centers and $C\not\subset f(E_f)$,
while \eqref{E363} requires some explanation.

Let $\Bl_{\wt{P}} \wt{X}$ and $\Bl_P X$ be the blowups of
$\wt{X}$ and $X$ along $\wt{P}$ and $P$, respectively.
We have the commutative diagram
\begin{equation}\label{E364}
\begin{CD}
\Bl_{\wt{P}} \wt{X} @>g>> \Bl_P X\\
@V\wt{\pi}VV @V{\pi}VV\\
\wt{X} @>f>> X
\end{CD}
\end{equation}
Let $\wt{\pi}_*^{-1} \wt{D}\subset \Bl_{\wt{P}} \wt{X}$ and
$\pi_*^{-1} D\subset \Bl_P X$ be the proper
transforms of $\wt{D}$ and $D$ under $\wt{\pi}$ and $\pi$,
respectively. Obviously, we have the commutative diagram
\begin{equation}\label{E365}
\begin{CD}
\wt{\pi}_*^{-1} \wt{D} @>g>> \pi_*^{-1} D\\
@V\wt{\pi}VV @V{\pi}VV\\
\wt{D} @>f>> D
\end{CD}
\end{equation}
and hence $g(\wt{\pi}_*^{-1} \wt{D}) = \pi_*^{-1} D$. Consequently,
\begin{equation}\label{E366}
g\left(\wt{C}\cap \wt{\pi}_*^{-1} \wt{D}\right)
\subset C\cap \pi_*^{-1} D
\end{equation}
and \eqref{E363} follows.

We conclude from \eqref{E361}-\eqref{E363} that
for every $C\subset
X$ satisfying \eqref{E351} and $C\not\subset f(E_f)$,
\begin{equation}\label{E353}
2 g(\wt{C}) - 2 + i_\wt{X}(\wt{C}, \wt{D}, \wt{P}) < (K_{\wt{X}} +
\wt{D}) \wt{C}.
\end{equation}
So the conclusion of the lemma holds for $(X, D, P)$ as long as it
holds for $(\wt{X}, \wt{D}, \wt{P})$, where $\wt{P} \subset
\wt{D}\cap \wt{F}$ and $\supp(\wt{F})$ has simple normal crossing.
Therefore, we may simply assume that $F$ has
simple-normal-crossing support at the very beginning.

Let $X'$ be the blowup of $X$ along $P$ and $\wt{X}$ be a
desingularization of $X'$. We claim that \eqref{E352} holds for
$f: \wt{X}\to X$. This follows from the following explicit
construction of $\wt{X}$.

Let $P = \cup_{i = 1}^\alpha \mu_i P_i$, where $P_i$'s are the
irreducible components of $P$ and $\mu_i$ is the multiplicity of
$P_i$ in $P$.

If $\mu_i = 1$ for each $i$, i.e., $P$ is reduced, $\wt{X}$
can be constructed by subsequently blowing up $X$ along
$P_1, P_2, ..., P_\alpha$. That is, we first blow up $X$ along
$P_1$ to obtain $X_1$, then blow up $X_1$ along the proper
transform of $P_2$ to obtain $X_2$, next blow up $X_2$ along the
proper transform of $P_3$ to obtain $X_3$ and so on. We obtain a
sequence of blowups:
\begin{equation}\label{E354}
\wt{X} = X_\alpha \xrightarrow{f_\alpha}
X_{\alpha-1}\xrightarrow{f_{\alpha-1}} ...\xrightarrow{f_3}
X_2\xrightarrow{f_2} X_1\xrightarrow{f_1} X_0 = X.
\end{equation}
It is not hard to see that \eqref{E354} is a sequence of blowups
with smooth centers due to our assumption on $F$ and $D$.
Therefore, $\wt{X}$ is smooth and it is easy to check
\eqref{E352} by observing that
\begin{equation}\label{E355}
K_{X_i} + D^{(i)} = f_i^*(K_{X_{i-1}} + D^{(i-1)})
\end{equation}
for each $i$, where $D^{(i)}$ is the
proper transform of $D$ under the map $X_i\to X$.

The construction of $\wt{X}$ becomes more complicated
if $P$ is nonreduced. Note that at a general point $p\in P_i$, $P$ is
locally given by $x = y^{\mu_i} = 0$, where $x = 0$ defines $D$
and $y^\mu = 0$ defines $F$ for some $\mu \ge
\mu_i$. So if $\mu_i > 1$, $X'$ has a singularity $p'$ over
$p$, which is locally given by $xz = y^{\mu_i}$. This suggests
that $\wt{X}$ is locally constructed over $p$ by blowing up
$X$ along $P_i$ $\mu_i$ times (see \remref{REM007} for
explanations). Correspondingly, the blowup $f_i:
X_i \to X_{i-1}$ in \eqref{E354} is replaced by a sequence of
$\mu_i$ blowups:
\begin{equation}\label{E356}
X_i = X_{i, \mu_i} \xrightarrow{f_{i,\mu_i}} X_{i, \mu_i - 1}
\xrightarrow{f_{i,{\mu_i - 1}}}
... \xrightarrow{f_{i, 1}} X_{i, 0} = X_{i-1}
\end{equation}
where $f_{i,1}: X_{i, 1}\to X_{i-1}$ is the blowup of $X_{i-1}$
along $\wt{P}_i$ and $f_{i, j}: X_{i, j} \to X_{i, j-1}$ is the
blowup of $X_{i, j-1}$ along $E_{i, j-1} \cap D^{(i, j-1)}$ for
$j > 1$. Here $\wt{P}_i$ is the proper transform of $P_i$ under
the map $X_{i-1}\to X$, $E_{i, j-1}$ is
the exceptional divisor of $f_{i, j-1}$ and $D^{(i, j-1)}$
is the proper transform of $D$ under the map $X_{i, j-1}\to
X$. Obviously, \eqref{E356} is a sequence of blowups with smooth
centers and it is easy to check that
\begin{equation}\label{E357}
K_{X_{i, j}} + D^{(i, j)} = f_{i,j}^*(K_{X_{i, j-1}} + D^{(i, j-1)})
\end{equation}
for each $j$. And hence \eqref{E355} and \eqref{E352} follow.

It follows from our construction of $f: \wt{X} \to X$ that
$f^{-1}(P)$ is a Cartier divisor of $\wt{X}$.
Hence $f$ factors through $X'$:
\begin{equation}\label{E358}
f: \wt{X} \xrightarrow{g} X' \to X
\end{equation}
by the universal property of blowups. Therefore, $g: \wt{X}\to X'$ is a
desingularization of $X'$. Let $\wt{D}$ and $D'$ be
the proper transforms of $D$ under the maps $\wt{X}\to X$ and $X'\to X$,
respectively. Then
\begin{equation}\label{E367}
i_\wt{X}(\wt{C}, \wt{D}) \le i_{X'}(C',
D') = i_X(C, D, P),
\end{equation}
where $\wt{C}$ and $C'$ are the proper
transforms of $C$ under the maps $\wt{X}\to X$ and $X'\to X$,
respectively, and we assume that $C\subset X$ is a curve not
contained in the image of the exceptional locus of $f$. Therefore,
\begin{equation}\label{E359}
\begin{split}
&\quad
2 g(\wt{C}) - 2 + i_\wt{X}(\wt{C}, \wt{D}) - (K_{\wt{X}} +
\wt{D})\wt{C}\\
&\le 2 g(C) - 2 + i_X(C, D, P) - (K_X + D)C.
\end{split}
\end{equation}
Therefore, to prove the lemma for $(X, D, P)$, it suffices to prove
it for $(\wt{X}, \wt{D}, \emptyset)$.
This reduces the lemma to the case $P = \emptyset$. Note that we can
no longer assume $D_k$ to be a general member of a BPF linear system
but we do have that $D$ is of simple normal crossing and that is all we
need in order to apply \lemref{LEM003}.
\end{proof}

\begin{rem}\label{REM006}
If $P$ is reduced, $g:\wt{X}\to X'$ is a small morphism, i.e., the
exceptional locus of $g$ has codimension at least two in $\wt{X}$.
Note that we may blow up $X$ along $P_i$'s in an arbitrary order
of $\{P_i\}$. By choosing a different order of $\{P_i\}$, we usually
arrive at a different desingularization of $X$, which is a
flop of $\wt{X}$. For example, let $X$ be a threefold and $P = P_1\cup
P_2$ be a curve in $X$. Suppose that $P_1$ and $P_2$ meet at a
point $p$ which is a node of the curve $P$. Then $X'$ has a
rational double point $p'$ over $p$. Let $X_{12}$ be the blowup
of $X$ along $P_1$ followed by the blowup along $P_2$ and
$X_{21}$ be the blowup of $X$ along $P_2$ followed by the blowup
along $P_1$. Then the rational map $X_{12}\to X_{21}$ induced by
\begin{equation}\label{E327}
X_{12}\xrightarrow{} X' \xleftarrow{} X_{21}
\end{equation}
is a flop.
\end{rem}

\begin{rem}\label{REM007}
To resolve the singularity $X = \{ xz = y^\mu \}$, we first blow up
$X$ along $P = \{ x = y = 0\}$ to obtain $X'$, which
has a singularity given by $x' z = y^{\mu-1}$ on $P' = \{x'
= y = 0\}$ with $x' = x/y$.
Blow up $X'$ along $P'$ to obtain $X''$, which has a singularity given
by $x'' z = y^{\mu-2}$ with $x'' = x'/y$. Continue this process and we
resolve the singularity $xz = y^{\mu}$ by blowing up $X$ along $P$
$\mu$ times.
\end{rem}

\subsection{Proof of \thmref{THM002}}\label{SEC002003}

If $C\subset D$, say $C\subset D_\alpha$ for some $\alpha\in I$,
then we may replace $(X, D)$ by $(X', D')$ with $X' = D_\alpha$ and
$D' = (D-D_\alpha)|_{D_\alpha}$ while observing that $\{D_i: i\ne \alpha\}$
satisfies $\C_{X'}(\varphi, 0, n-2)$ on $X'$. So let us assume that
$C$ meets $D$ properly. It suffices to prove the following statement.

\begin{prop}\label{PROP007}
Under the hypotheses of \thmref{THM002}, the following holds:
for each $J\subset I$ with $|J| = n - r$, the reduced
irreducible curves $C\not\subset D\cup F$ satisfying 
\begin{equation}\label{E370}
2 g(C) - 2 + i_X(C, \sum_{i\not\in J} D_i, P) < \varphi(C)
\end{equation}
are at most $(n - r)$-filling for $1\le r\le n$.
\end{prop}

\begin{proof}
We argue by induction on $r$.

By \lemref{LEM001}, for each subset $J\subset I$
with $|J| = n - 1$, the curves $C$ satisfying
\begin{equation}\label{E389}
\left(K_X + \sum_{i\not\in J} D_i\right) C < \varphi(C)
\end{equation}
are at most $(n-1)$-filling; and by \lemref{LEM004}, the reduced
irreducible curves $C$ satisfying
\begin{equation}\label{E388}
2 g(C) - 2 + i_X(C, \sum_{i\not\in J} D_i, P) 
< \left(K_X + \sum_{i\not\in J} D_i\right) C
\end{equation}
are at most $(n-1)$-filling. Therefore, the reduced irreducible
curves $C$ satisfying \eqref{E370} are at most
$(n-1)$-filling for $|J| = n - 1$. So the proposition holds for $r
= 1$.

We use the notation $\overline{\sum_{i\not\in J} D_i}$ to denote the
closure of the locally noetherian subscheme of $X$ swept out by the
reduced irreducible curves $C\not\subset D\cup F$
satisfying \eqref{E370}. The statement of the proposition is
equivalent to that
\begin{equation}\label{E383}
\dim \overline{\sum_{i\not\in J} D_i} \le n - r
\end{equation}
for each $J\subset I$ with $|J| = n - r$.

Suppose that
\begin{equation}\label{E411}
\dim \overline{\sum_{i\not\in J} D_i} \le n - r + 1
\end{equation}
for each $J\subset I$ with $|J| = n - r + 1$. We want to show that
\eqref{E383} holds. Assume the contrary:
\begin{equation}\label{E384}
\dim \overline{\sum_{i\not\in J} D_i} = n - r + 1
\end{equation}
for some $J\subset I$ with $|J| = n - r$.
Fix $k\not\in J$ and let $J' = J\cup \{k\}$. Then according
to the inductive hypothesis,
\begin{equation}\label{E385}
\dim \overline{\sum_{i\not\in J'} D_i} \le n - r + 1.
\end{equation}
And since
\begin{equation}\label{E371}
\overline{D_k + \sum_{i\not\in J'} D_i} =
\overline{\sum_{i\not\in J} D_i}\subset \overline{\sum_{i\not\in J'} D_i},
\end{equation}
\begin{equation}\label{E372}
\dim \overline{D_k + \sum_{i\not\in J'} D_i}
= \dim \overline{\sum_{i\not\in J'} D_i}.
\end{equation}
Therefore, the irreducible components of $\overline{D_k + \sum_{i\not\in
J'} D_i}$ of dimension $n - r + 1$ must also be the
components of $\overline{\sum_{i\not\in J'} D_i}$. Now if we
choose $D_k' \in \P \L_k$ and $D_k' \ne D_k$, we still have
\begin{equation}\label{E373}
\overline{D_k' + \sum_{i\not\in J'} D_i}\subset
\overline{\sum_{i\not\in J'} D_i}
\end{equation}
and
\begin{equation}\label{E374}
\dim \overline{D_k' + \sum_{i\not\in J'} D_i}
= \dim \overline{\sum_{i\not\in J'} D_i}.
\end{equation}
Therefore, as $D_k$ varies in the linear series $\P \L_k$, the
$(n-r+1)$-dimensional irreducible components of $\overline{D_k +
\sum_{i\not\in J'} D_i}$ vary as irreducible components of
$\overline{\sum_{i\not\in J'} D_i}$. When we fix $\sum_{i\not\in
J'} D_i$, the scheme $\overline{\sum_{i\not\in J'} D_i}$ is
fixed and hence so are its components. Consequently, the
$(n-r+1)$-dimensional components of $\overline{D_k +
\sum_{i\not\in J'} D_i}$ remain fixed as we vary $D_k$
and fix the rest of $D_i$'s. Therefore,  the
$(n-r+1)$-dimensional components of $\overline{D_k +
\sum_{i\not\in J'} D_i}$ and $\overline{D_k' +
\sum_{i\not\in J'} D_i}$ are the same. From now on, we will
drop all the components of dimension less than $n - r + 1$ from
the definition of $\overline{\sum_{i\not\in J} D_i}$, i.e.,
we will pretend that $\overline{\sum_{i\not\in J} D_i}$ is pure of
dimension $n-r+1$. Then we have
\begin{equation}\label{E375}
\overline{D_k + \sum_{i\not\in J'} D_i} = \overline{D_k' +
\sum_{i\not\in J'} D_i}
\end{equation}
for any general choice of $D_k, D_k'\in \P \L_k$. Note that $k$ is
chosen arbitrarily from $I\backslash J$. So
\eqref{E375} actually implies that
\begin{equation}\label{E376}
\overline{\sum_{i\not\in J} D_i} = \overline{\sum_{i\not\in J} D_i'}
\end{equation}
for any general choices of $D_i, D_i'\in \P\L_i$. Or equivalently,
$\overline{\sum_{i\not\in J} D_i}$ remains fixed even if we vary all
$D_i$'s at the same time.

A word of warning is in order for a potential misunderstanding of
\eqref{E376}. It does not imply that a curve $C\subset X$
satisfying \eqref{E370} will also satisfy
\begin{equation}\label{E377}
2 g(C) - 2 + i_X(C, \sum_{i\not\in J} D_i', P) < \varphi(C).
\end{equation}
The curves satisfying \eqref{E370} and \eqref{E377} may very well
be different but they do sweep out the same subscheme in $X$, if we
assume that \eqref{E384} holds.

Let $Y$ be an irreducible component of $\overline{\sum_{i\not\in J}
D_i}$ and $f: \wt{Y}\to Y\subset X$ be a desingularization of
$Y$. Let $\wt{D}_i = f^* D_i$, $\wt{F} = f^* F$ and $\wt{P} =
f^{-1}(P)$. Since $Y$ is independent of the choices of $D_i$ by
\eqref{E376}, $\wt{D}_i$ is a very general member of a BPF
linear series on $\wt{Y}$.

For a reduced irreducible curve $C\subset Y$ and
$C\not\subset D$, we have $g(\wt{C}) = g(C)$
and
\begin{equation}\label{E382}
i_{\wt{Y}}(\wt{C}, \sum_{i\not\in J}\wt{D}_i, \wt{P})\le i_X(C,
\sum_{i\not\in J} D_i, P)
\end{equation}
by the same argument for \eqref{E363},
where $\wt{C}$ is the proper transform of $C$ under $f$.
Hence
\begin{equation}\label{E378}
2 g(\wt{C}) - 2 + i_{\wt{Y}}(\wt{C}, \sum_{i\not\in J}\wt{D}_i,
\wt{P})
\le 2 g(C) - 2 + i_{X}(C, \sum_{i\not\in J} D_i, P).
\end{equation}
Therefore, $\wt{Y}$ is covered by
the reduced irreducible curves $\wt{C}$ satisfying 
\begin{equation}\label{E379}
2 g(\wt{C}) - 2 + i_{\wt{Y}}(\wt{C}, \sum_{i\not\in J}\wt{D}_i,
\wt{P}) < \varphi(C).
\end{equation}
Since such curves $\wt{C}$ cover $\wt{Y}$, we must have
\begin{equation}\label{E380}
(K_\wt{Y} + \sum_{i\not\in J}\wt{D}_i)\wt{C} \ge \varphi(C)
\end{equation}
by \lemref{LEM001}. Therefore,
\begin{equation}\label{E381}
2 g(\wt{C}) - 2 + i_{\wt{Y}}(\wt{C}, \sum_{i\not\in J}\wt{D}_i,
\wt{P}) < (K_\wt{Y} + \sum_{i\not\in J}\wt{D}_i)\wt{C}.
\end{equation}
This contradicts with \lemref{LEM004}, which says the curves satisfying
\eqref{E381} are at most $(\dim\wt{Y} - 1)$-filling in
$\wt{Y}$.
\end{proof}

\section{A Special Case of \thmref{THM003}}\label{SEC003}

Due to the technicality of the proof of \thmref{THM003}, we
feel it is better to first work out a special case. So we will delay
the proof of \thmref{THM003} to the next section. Here instead we will
study a special case, as the major ingredients of the proof are
already present in this special case.

\subsection{The complement of a surface of degree 7 in $\P^3$}
\label{SEC003001}

Our proof of \thmref{THM003} consists of two main steps of degeneration.
\begin{description}
\item[Step 1]
First we specialize $D$ by degenerating $D_{I_\beta}$ to
$\sum_{i\in I_\beta} D_i$ for each $\beta\in\CB$.
In this way, we reduce the theorem to the case
$|I_\beta| = 1$ for all $\beta\in\CB$.
\item[Step 2]
Next we specialize $S$ by degenerating $D_{I_\alpha}$ to
$\sum_{i\in I_\alpha} D_i$ for each $\alpha\in \CA$.
Eventually, we reduce the theorem to
\thmref{THM002}, which we have already proved.
\end{description}

The second step of our proof relies on a concept developed in
\cite{C2}, called {\it virtual genus\/} of a curve lying on a
reducible variety. 
Here we will work out an example where $\CA = \emptyset$. So it only
involves the first step of degeneration.
The case $\dim X = 2$ was handled in \cite{C2}.
Let us consider the first nontrivial case in $X = \P^3$:
the complement of a very general surface $D$ of degree $7$.
We want to prove that
\begin{equation}\label{E202}
2 g(C) - 2 + i_{\P^3}(C, D) \ge \deg C
\end{equation}
for all reduced curves $C\subset \P^3$. Our inductive hypothesis is that
the theorem holds in $\dim X - |\CA| < 3$. In particular, we assume:

\begin{hyp}\label{HYP002}
On a very general sextic surface $S\subset\P^3$,
\begin{equation}\label{E201}
2 g(C) - 2 \ge \deg C
\end{equation}
for all reduced curves $C\subset S$.
\end{hyp}

As planned, the proof of \eqref{E202} is carried out by degenerating
$D$ to a union of seven planes in general position for which
\thmref{THM002} can be applied. However, due to technical
difficulties, we find it better to degenerate one plane at a time
instead of degenerating $D$ directly to seven planes. Namely, we will
first degenerate $D$ to a union of a sextic surface and a plane, then
to a union of a quintic surface and two planes and so on. The
argument for each step of degeneration is similar. We will illustrate
it by carrying out the first step of degeneration. That is, we will
prove \eqref{E202} under the inductive hypothesis that

\begin{hyp}\label{HYP003}
\eqref{E202} holds for $D = S\cup G\subset \P^3$
a very general union of a sextic surface $S$ and a plane $G$, i.e.,
\begin{equation}\label{E203}
2 g(C) - 2 + i_{\P^3}(C, S\cup G) \ge \deg C
\end{equation}
for all reduced curves $C\subset\P^3$.
\end{hyp}

Thus we will degenerate $D$ to $S\cup G$. The basic set up is as follows.

Let $Z = \P^3\times \Delta$ and $W \subset Z$ be a pencil of
surfaces in $\P^3$ of degree $7$ such that $W$ is irreducible
and $W_0 = S\cup G$.

Let $Y$ be a reduced flat family of curves over $\Delta$ with the
commutative diagram \eqref{E204}. Our goal is to prove that
\begin{equation}\label{E205}
2 g(Y_t) - 2 + i_Z(\pi_* Y_t, W) \ge \deg \pi_* Y_t
\end{equation}
for $t\ne 0$ by analyzing what happens on the central fiber.
If $\pi_* Y_0$ is reduced and meets $W_0$ properly, then
\begin{equation}\label{E206}
2 g(Y_0) - 2 + i_{\P^3}(\pi_* Y_0, W_0) \ge \deg \pi_* Y_0
\end{equation}
by \Hypref{HYP003} and \eqref{E205} follows due to the
obvious semi-continuity of its LHS. However, $\pi_* Y_0$ could very
well be nonreduced and even worse, it may fail to meet $W_0$
properly. Overcoming these difficulties is essentially what our proof
of \thmref{THM003} is about.

It turns out that the biggest problem we could have is that $\pi_* Y_0$
has a component contained in $G$ while the other problems mentioned
above can be easily overcome. Fortunately, there is a very common
construction we can use to resolve the problem.

\subsection{Construction of a fan}\label{SEC003002}

Let $f: \wt{Z}\to Z$ be the blowup of $Z$ along $G$ and let $R$ be the
exceptional divisor (see \figref{FIG001}).
Obviously, the central fiber $\wt{Z}_0$ of
$\wt{Z}$ is the union $R\cup X$, where $R$ is a $\P^1$ bundle
over $G$ and $R\cap X = G$. Actually, $R$ is the
projectivization of the normal bundle $N_{G/Z}$ of $G$ in
$Z$. It is not hard to see that
\begin{equation}\label{E207}
N_{G/Z} = \CO_G \oplus \CO_G(G) = \CO_{\P^2} \oplus \CO_{\P^2}(1).
\end{equation}
So $R$ is the $\P^1$ bundle over $\P^2$ given by
$\P(\CO_{\P^2}\oplus \CO_{\P^2}(1))$.

Such construction has been extensively used by Z. Ran in his study of
Severi varieties of plane curves \cite{R} and more recently by C. Ciliberto
and R. Miranda in their works on Nagata conjecture \cite{CM1} and
\cite{CM2}. In their cases, $Z = \P^2\times \Delta$ and $G\subset
Z_0\isom \P^2$ is a line. Then blow up $Z$ along $G$ and we will
obtain the ruled surface $R\isom \F_1$ over $G$ on the central fiber
$\wt{Z}_0$ as the exceptional divisor.
Such construction can be easily generalized. More generally, $Z$ can
be an arbitrary flat family of varieties over $\Delta$ and
$G\subset Z_0$ be a closed subscheme of the central fiber
$Z_0$. We blow up $Z$ along $G$ to obtain $\wt{Z}$. The
central fiber $\wt{Z}_0$ consists of the proper transform of $Z_0$
and the exceptional divisor $R$, which is a $\P^n$ bundle over
$G$. We will use the terminology of Ran to call $\wt{Z}_0$
a {\it fan\/}. 

The main purpose of constructing a fan $\wt{Z}_0$ is to study the
infinitesimal behavior of a flat family $Y\subset Z$ in the
neighborhood of $G$. In our case, we want to ``separate'' $\pi(Y)$
and $W$ along $G$ via the construction of $\wt{Z}$.

Let $\wt{W}$ be the proper transform of $W$ under the blowup
$\wt{Z} \to Z$. The central fiber $\wt{W}_0$ of $\wt{W}$
is the union of the proper transform of $S$, which we still denote
by $S$, and a surface $\wt{G}\subset R$ (see \figref{FIG001}).

\psset{unit=0.01in,linewidth=1.2pt,dash=3pt 2pt,dotsep=2pt}


\begin{figure}[ht]
\centering
\begin{pspicture}(-170,-300)(170,10)

\psline(-50,0)(150,0)
\psline(-150,-100)(50,-100)
\psline(-50,0)(-150,-100)
\psline(150,0)(50,-100)
\psline(-100,-50)(100,-50)

\uput[225](-150,-30){\small $Z_0$}
\uput[225](-150,-100){\small $X$}

\uput[270](65,-50){\small $G$}

\parametricplot{0}{360}{50 t cos mul 35 t sin mul add 
-50 35 t sin mul add}

\uput[225](-60,-75){\small $S$}

\uput[225](-60,-255){\small $S$}

\psline[doubleline=true,doublesep=1.5pt,arrowsize=2pt
2]{->}(0,-105)(0,-135)

\uput[45](100,-140){\small $R$}
\uput[225](-150,-160){\small $\wt{Z}_0$}

\psline(-100,-140)(100,-140)
\psline(-100,-230)(100,-230)
\psline[linestyle=dashed](-50,-180)(100,-180)
\psline(100,-180)(150,-180)
\psline(-150,-280)(50,-280)

\uput[225](-150,-280){\small $X$}

\uput[270](65,-230){\small $G$}

\parametricplot[linestyle=dashed]{0}{180}{50 t cos mul 35 t sin mul add 
-230 35 t sin mul add}

\parametricplot{180}{360}{50 t cos mul 35 t sin mul add 
-230 35 t sin mul add}

\parametricplot{0}{180}{50 t cos mul -230 60 t sin mul add}

\uput[90](0,-170){\small $\wt{G}$}

\psline[linestyle=dashed](-50,-180)(-100,-230)
\psline(-100,-230)(-150,-280)

\psline(-100,-140)(-100,-230)
\psline(100,-140)(100,-230)

\psline(150,-180)(50,-280)

\end{pspicture}

\caption{The blowup of \(Z = \P^3\times\Delta\) along a plane
$G\subset Z_0$}\label{FIG001}
\end{figure}
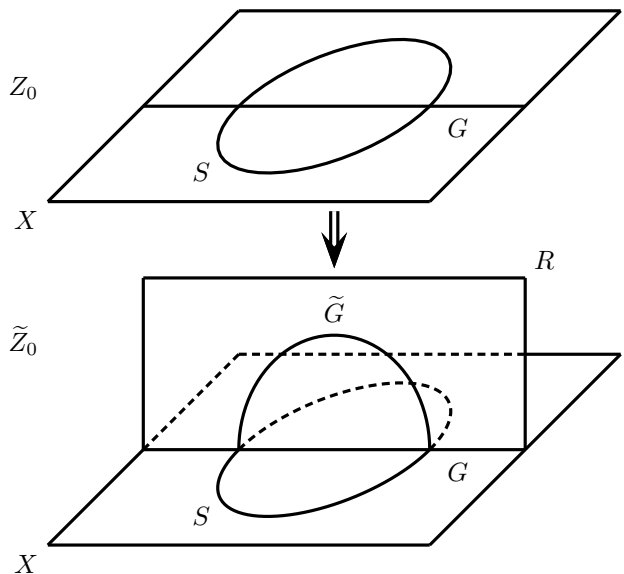

It is not hard to figure
out what $\wt{G}$ is. First of all, we obviously have
\begin{equation}\label{E208}
S\cap G = \wt{G}\cap G.
\end{equation}
Second, it is not hard to see that $f_*(\wt{G}) = G$, where $f:
R\to G$ is the projection. Indeed, $\wt{G}$ is a member of the
linear series
\begin{equation}\label{E209}
\P H^0(\CO_R(G)\tensor f^*(\CO_G(7)))
\end{equation}
where $\CO_R(G)$ is the tautological line bundle of $R$ and $\CO_G(1)$
is the hyperplane bundle of $G\isom \P^2$.
In addition, for a generic choice of the pencil $W$, the
corresponding $\wt{G}$ is a general member of the linear series
\eqref{E209}.

Now we have a rational map $Y\to \wt{Z}$ factoring through $Z$. After
resolving the indeterminacies of this map, we obtain
\begin{equation}\label{E210}
\begin{CD}
\wt{Y} @>\wt{\pi}>> \wt{Z}\\
@VVV @VVfV\\
Y @>\pi>> Z.
\end{CD}
\end{equation}
Since the blowup $f$ did nothing to the general fibers, we have
\begin{equation}\label{E211}
2 g(Y_t) - 2 + i_Z(\pi_* Y_t, W) =
2 g(\wt{Y}_t) - 2 + i_\wt{Z}(\wt{\pi}_* \wt{Y}_t, \wt{W})
\end{equation}
for $t\ne 0$.

We are trying to bound the RHS of \eqref{E211} with the
information on $\wt{Y_0}$, which is a curve lying on the union of two
smooth varieties $X$ and $R$ meeting transversely. There is a result
\cite[Theorem 1.17]{C2} dealing with this situation. 
We will put it in a more general form suitable for our purpose.

\begin{prop}\label{PROP002}
Let $X$ be a flat family of projective varieties over
$\Delta$, whose general fibers $X_t$ are irreducible and smooth
and whose central fiber
$X_0 = D = \cup_{i\in I} D_i$ is of normal crossing 
along $\partial D_i$ for each $i\in I$, where
$D_i$ are irreducible components of $X_0$ and we write
\begin{enumerate}
\item $\partial D_i = D_i\cap (\cup_{j\ne i} D_j)$ for $i\in I$,
\item more generally, $D_J = \cap_{j\in J} D_j$,
$\partial D_J = D_J\cap (\cup_{i\not\in J} D_i)$
for every $J\subset I$,
\item and $\partial D = \cup_{i\in I} \partial D_i = \cup_{|J| = 2} D_J$.
\end{enumerate}

Let $Q = \partial D \cap X_{sing}$
be the singular locus of $X$ along $\partial D$ and
suppose that $\dim (Q\cap D_J)\le \dim D_J - 1$ for every $J\subset I$.

Let $W \subset X$ be an effective divisor of $X$ that is flat over
$\Delta$. Suppose that $W$ meets $D_J$ and
$D_J\cap Q$ properly in $X$ for each $J\subset I$.

Let $Y$ be a reduced flat family of curves over $\Delta$ with the
commutative diagram:
\begin{equation}\label{E212}
\begin{CD}
Y @>\pi>> X\\
@VVV @VVV\\
\Delta @>>> \Delta
\end{CD}
\end{equation}
where $\pi: Y\to X$ is a proper map and $\Delta\to \Delta$ is a base
change.

Then
\begin{equation}\label{E213}
2 g(Y_t) - 2 + i_X(\pi_* Y_t, W)\ge
\sum_{\Gamma\subset \pi_* Y_0} \mu_\Gamma \Phi_{X,W}(\Gamma)
\end{equation}
where we sum over all irreducible components
$\Gamma\subset \pi_* Y_0$,
$\mu_\Gamma$ is the multiplicity of $\Gamma$ in $\pi_* Y_0$ and
\begin{equation}\label{E214}
\Phi_{X,W}(\Gamma) =
2 g(\Gamma) - 2 + i_X(\Gamma, D\cup W, \partial D_J\cap Q)
\end{equation} 
with $J = \{j\in I: \Gamma\subset D_j\}$.

Furthermore, suppose that there exists a morphism $f: X\to \hat X$
which blows down $D_\alpha$ for some $\alpha\in I$ while inducing
isomorphisms on $\cup_{i\ne\alpha} D_i$ and $X_t$ for $t\ne 0$.
Let $\hat W = f_* W$. Then
\begin{equation}\label{E215}
\begin{split}
&\quad \quad \quad 2 g(Y_t) - 2 + i_X(\pi_* Y_t, W)
\\&\ge 
\sum_{\genfrac{}{}{0pt}{}{\Gamma\subset D_\alpha}{f_* \Gamma\ne 0}}
\mu_\Gamma \Phi_{X,W}(\Gamma)
+ \sum_{\Gamma\not\subset D_\alpha} \mu_\Gamma \Phi_{\hat X,\hat W}(f_* \Gamma)
\end{split}
\end{equation}
where we sum over the irreducible components $\Gamma\subset\pi_* Y_0$
not contracted by $f$.

\noindent {\bf Modification:}
if we assume that $Y$ is irreducible,
then \eqref{E213} and \eqref{E215} can be slightly improved as
follows. Let $U_Y$ be the union of the irreducible components $\Gamma\subset
\pi_* Y_0$ with the following properties:
\begin{enumerate}
\item $\Phi_{X,W}(\Gamma) < 0$;
\item $\Gamma\cap \partial D = \emptyset$;
\item $B \cdot \Gamma\ge 0$ for every irreducible component $B$
of $W$ satisfying that $\Gamma\subset B$ and $\pi(Y)\not\subset B$.
\end{enumerate}
If $U_Y\subsetneq \supp(\pi_* Y_0)$, we may exclude $\Gamma\subset U_Y$ in
the summations of \eqref{E213} and \eqref{E215}, i.e., 
\eqref{E213} and \eqref{E215} continue to
hold when we sum over $\Gamma\not\subset U_Y$.
\end{prop}

\begin{rem}\label{REM001}
If we take $W = \emptyset$, \eqref{E213} becomes
\begin{equation}\label{E216}
g(Y_t) \ge g_Q^{vir}(\pi_* Y_0)
\end{equation}
with $g_Q^{vir}(\pi_* Y_0)$ the {\it virtual genus\/} of $\pi_* Y_0$
with respect to $Q$ defined in \cite{C2} and this is
exactly what \cite[Theorem 1.17]{C2} says.
\end{rem}

We will prove \propref{PROP002} much later in \ssecref{SEC004003}.
Now let us apply it to our situation with $\wt{Y}$ being the family of
curves and $f:\wt{Z}\to Z$ being the map blowing down $R$.
We obtain
\begin{equation}\label{E217}
\begin{split}
&\quad\quad\quad 2 g(Y_t) - 2 + i_Z(\pi_* Y_t, W)\\
&\ge \sum_{\genfrac{}{}{0pt}{}{\Gamma\subset R}{f_* \Gamma\ne 0}}
\mu_\Gamma\Phi_{\wt{Z},\wt{W}}(\Gamma)
+ \sum_{\Gamma\not\subset R} \mu_\Gamma \Phi_{Z,W}(\Gamma)
\end{split}
\end{equation}
where we sum over the irreducible components $\Gamma\subset \wt{\pi}_*
\wt{Y}_0$ not contracted by $f$.

By the definition of $\Phi$, we have
\begin{align}\label{E227}
\Phi_{\wt{Z},\wt{W}}(\Gamma) &=
2 g(\Gamma) - 2 + i_R(\Gamma, G \cup \wt{G}) &&\text{ for }
\Gamma\subset R,\\
\label{E400}
\Phi_{Z,W}(\Gamma) &= 2 g(\Gamma) - 2 &&\text{ for } \Gamma\subset S,
\text{ and}\\
\label{E009}
\Phi_{Z,W}(\Gamma) &= 2 g(\Gamma) - 2 + i_X (\Gamma, S\cup G)
&& \text{ for } \Gamma\not\subset S\cup R.
\end{align}

The RHS's of \eqref{E400} and \eqref{E009} can be easily
bounded by the inductive hypothesis, i.e.,
\begin{equation}\label{E218}
2 g(\Gamma) - 2 \ge \deg \Gamma
\end{equation}
for $\Gamma\subset S$ by \eqref{E201} and
\begin{equation}\label{E390}
2 g(\Gamma) - 2 + i_X (\Gamma, S\cup G) \ge \deg \Gamma
\end{equation}
for $\Gamma\subset X$ by \eqref{E203}.
Therefore, in order to prove \eqref{E205}, it suffices to show
\begin{equation}\label{E220}
2 g(\Gamma) - 2 + i_R(\Gamma, G \cup \wt{G}) \ge \deg f_*\Gamma
\end{equation}
for $\Gamma\subset R$ and $f_*\Gamma\ne 0$,
since $\deg \pi_* Y_t = \sum \mu_\Gamma \deg f_*\Gamma$.
So we have eventually reduced \eqref{E202},
which is a statement on the log
pair $(X, D)$, to \eqref{E220}, which is a statement on the log
pair $(R, G\cup \wt{G})$. Note that $R$ is a $\P^1$ bundle over
$G\isom \P^2$. To prove \eqref{E220}, we need to further
degenerate $\wt{G}$. This will lead to the construction of a fan by
blowing up $R\times \Delta$ along $F\subset R$, where $F = f^* L$ is the
pullback of a line $L\subset G$ under the projection $f: R\to G$. The
exceptional divisor of this blowup is a $\P^1$ bundle $R'$ over $F$
and $R'$ can also be regarded as a $\P^1\times \P^1$ bundle over
$L\isom \P^1$ (see \figref{FIG002}). 
We call $R'$ a {\it projective tower\/} over $L$.


\begin{figure}[ht]
\centering
\begin{pspicture}(-230,-350)(230,50)

\psline(-120,0)(-20,0)
\psline(-200,-80)(-100,-80)
\psline(-120,0)(-200,-80)
\psline(-20,0)(-100,-80)
\psline(-70,0)(-150,-80)

\uput[225](-200,-80){\small $X$}
\uput[330](-110,-40){\small $G$}
\uput[225](20,-80){\small $X$}
\uput[330](110,-40){\small $G$}
\uput[45](150,35){\small $R$}

\rput[lt]{0}(-200,-100){
\parbox{1.7in}{\small Blow up $X\times\Delta$ along $G\isom\P^2$ lying
on the central fiber}
}

\psline[doubleline=true,doublesep=1.5pt,arrowsize=2pt
2]{->}(-15,-40)(15,-40)

\psline(100,0)(115,0)
\psline[linestyle=dashed](115,0)(150,0)
\psline(150,0)(200,0)

\psline(20,-80)(120,-80)
\psline(100,0)(20,-80)
\psline(200,0)(120,-80)
\psline(150,0)(70,-80)

\psline(70,-80)(70,-45)
\psline(150,0)(150,35)
\psline(70,-45)(150,35)

\rput[lt]{0}(20,-100){
\parbox{1.7in}{\small We obtain the exceptional divisor $R\isom\P_G\E$
where $\E = \CO \oplus \CO(1)$.}
}

\rput[lt]{0}(-230,-150){
\begin{pspicture}(-230,-150)(230,50)
\psline(-120,0)(-20,0)
\psline(-200,-80)(-100,-80)
\psline(-120,0)(-200,-80)
\psline(-20,0)(-100,-80)
\psline(-70,0)(-150,-80)

\rput[lt]{0}(-200,-100){
\parbox{1.7in}{\small Blow up $R\times\Delta$ along $F\isom\P_L \E$ lying
on the central fiber}
}

\psline[doubleline=true,doublesep=1.5pt,arrowsize=2pt
2]{->}(-15,-40)(15,-40)

\psline(100,0)(115,0)
\psline[linestyle=dashed](115,0)(150,0)
\psline(150,0)(200,0)

\psline(20,-80)(120,-80)
\psline(100,0)(20,-80)
\psline(200,0)(120,-80)
\psline(150,0)(70,-80)

\psline(70,-80)(70,-45)
\psline(150,0)(150,35)
\psline(70,-45)(150,35)

\uput[225](-200,-80){\small $R$}
\uput[330](-110,-40){\small $F$}
\uput[225](20,-80){\small $R$}
\uput[330](110,-40){\small $F$}
\uput[45](150,35){\small $R'$}

\rput[lt]{0}(20,-100){
\parbox{1.7in}{\small We obtain the exceptional divisor
$R'\isom\P_L\E\times_L \P_L\E$.}
}
\end{pspicture}
}

\end{pspicture}

\caption{The projective tower
\(R' = \P_L(\E, \E)\)
over \(L = \P^1\)
}
\label{FIG002}
\end{figure}

\subsection{Projective Tower}\label{SEC003003}

Let $X$ be a scheme and $\E_1, \E_2, ..., \E_n$ be vector bundles
over $X$. The {\it projective tower\/} of $\E_1, \E_2, ..., \E_n$
over $X$,
denoted by $\P_X(\E_1, \E_2, ..., \E_n)$ or
$\P(\E_1, \E_2, ..., \E_n)$ if $X$ is clear from the context,
is constructed inductively
as follows. First, let $X_0 = X$ and $X_1 = \P \E_1$. Next, let $X_2 =
\P (\pi_1^* \E_2)$ be the projectivization of the vector bundle
$\pi_1^* \E_2$ over $X_1$, where $\pi_1$ is the projection $X_1 \to X_0$.
Similarly, $X_3$ is the projectivization of the pullback
of $\E_3$ over $X_2$ and so on. Finally, $X_n$ is the
projectivization of the pullback of $\E_n$ over $X_{n-1}$ and we
call $X_n = \P(\E_1, \E_2, ..., \E_n)$ the projective tower of
$\E_1, \E_2, ..., \E_n$ over $X$. The order of $\E_1, \E_2, ...,
\E_n$ in the construction is not important, e.g., $\P(\E_1,
\E_2, ..., \E_n) = \P(\E_2, \E_1, ..., \E_n)$. Indeed, $\P(\E_1,
\E_2, ..., \E_n)$ is just the fiber product $\P\E_1 \times_X \P\E_2
\times_X ... \times_X \P \E_n$ over $X$. The Picard group of
$\P(\E_1, \E_2, ..., \E_n)$ is $\Pic X \oplus \BZ^n$, where
the $\BZ^n$ part is generated by $M_1, M_2, ..., M_n$ with $M_i$
the pullback of the tautological divisor of $\P \E_i$ under the
projection $\P(\E_1, \E_2, ..., \E_n)\to \P \E_i$.
We call $M_i$ the tautological divisors of $\P(\E_1, \E_2, ...,
\E_n)$.
Let $D\in \Pic X$ be a divisor on $X$ and $\pi$ be the projection
$\P(\E_1, \E_2, ..., \E_n)\to X$. Then the global sections of
a divisor $a_1 M_1 + a_2 M_2 + ... + a_n M_n + \pi^* D$ on
$\P(\E_1, \E_2, ..., \E_n)$ can be naturally identified with the
global sections of the vector bundle
\begin{equation}\label{E301}
\Sym^{a_1} \E_1^\vee \tensor \Sym^{a_2} \E_2^\vee \tensor ... \tensor
\Sym^{a_n} \E_n^\vee \tensor \CO_X(D)
\end{equation}
on $X$, where $\E_i^\vee$ is the dual of $\E_i$.

Before we proceed, we want to state a simple principle that was used throughout
\cite{C2} and will be used extensively here as well. We will glorify
it by calling it ``Riemann-Hurwitz on curves with marked
points''. Basically, it says the following.

\begin{prop}
\label{PROP003}
Let $\pi_\Gamma: \Gamma\to \Sigma$ be a surjective map between two
reduced curves $\Gamma$ and $\Sigma$
and let $M\subset \Sigma$
be a finite set of points on $\Sigma$. Let $\pi_\Gamma^{-1}(M)$ be
the (set-theoretical) inverse image of $M$ and we assume that
$\pi_\Gamma^{-1}(M)$ is contained in the nonsingular locus of
$\Gamma$. Then
\begin{equation}\label{E302}
2 g(\Gamma) - 2 + |\pi_\Gamma^{-1}(M)| \ge \gamma(2g(\Sigma) - 2 +
|M|)
\end{equation}
where $\gamma$ is the degree of $\pi_\Gamma$.
\end{prop}

This is more or less obvious by noting that the total ramification
index of the points in $\pi_\Gamma^{-1}(M)$ is at least $\gamma |M|
- |\pi_\Gamma^{-1}(M)|$.

\begin{prop}\label{PROP004}
Let $\E$ be the rank two vector bundle over $\P^k$ given by
$\E = \CO \oplus \CO(1)$ and 
\begin{equation}\label{E221}
R = \P(\underbrace{\E, \E, ..., \E}_{n-k}) = \P \E^{(n-k)}
\end{equation}
be the projective tower of $n-k$ $\E$'s over $\P^k$. Let $M_i$
be the pullback of the tautological divisor under the $i$-th
projection $p_i: R\to \P \E$ and let $\f$ be the pullback of the
hyperplane divisor of $\P^k$ under the projection $p: R\to
\P^k$. Suppose that $L$ is a very general member of the linear series
$|M + a \f|$ and $F_1, F_2, ..., F_b$ are
$b$ very general members of the linear series $|\f|$, where $M = M_1 +
M_2 + ... + M_{n-k}$ and $a\ge n-k\ge 0$. Then 
\begin{equation}\label{E222}
2 g(C) - 2 + i_R(C, L\cup M\cup F) \ge (a + b - 2n) (C\cdot \f)
\end{equation}
for every reduced irreducible curve $C\subset R$ satisfying $p_* C \ne
0$, where $F = F_1 + F_2 + ... + F_b$. Here we assume that $0 < k < 3$.
\end{prop}

In \eqref{E220}, $R = \P_{\P^2} \E$, $G\in \P
H^0(\CO_R(M))$ and $\wt{G}\in \P H^0(\CO_R(M + 7 \f))$. Therefore,
\eqref{E220} follows directly from the above
proposition by taking $n = 3$, $k = 2$, $a = 7$ and $b = 0$.

\begin{proof}[Proof of \propref{PROP004}]
First, let us prove the proposition when $a = n - k$. Let
$\Gamma =\supp(p_* C)$
be the reduced image of $C$ under the projection $p$ and
let $p_C: C\to \Gamma$ be the restriction of $p$ to $C$.
Applying \propref{PROP003} to the map $p_C$ by noticing
that
\begin{equation}\label{E223}
p_C^{-1}(\Gamma\cap F_\beta) \subset C\cap F_\beta, \text{ for }
\beta = 1,2,...,b,
\end{equation}
we obtain (see \figref{FIG003})
\begin{equation}\label{E224}
2 g(C) - 2 + i_R(C, F)
\ge (\deg p_C) \left(2g(\Gamma) - 2 + i_{\P^k}(\Gamma, p(F))\right).
\end{equation}
Obviously, $p(F) = p(F_1)\cup p(F_2)\cup ...\cup p(F_b)$ is the union
of $b$ hyperplanes of $\P^k$ in general position. Therefore,
\begin{equation}\label{E225}
2g(\Gamma) - 2 + i_{\P^k}(\Gamma, p(F)) \ge (b-2k) \deg \Gamma
\end{equation}
by \thmref{THM002}. Combining \eqref{E224} and \eqref{E225} yields
\begin{equation}\label{E226}
2 g(C) - 2 + i_R(C, F) \ge (b-2k) (C\cdot \f).
\end{equation}
And hence \eqref{E222} holds when $a = n -k$.


\begin{figure}[ht]
\centering
\begin{pspicture}(-230,-20)(230,100)

\psline(-200,0)(0,0)

\psline(-200,0)(-200,80)

\uput[90](-200,80){\small $F_1$}

\psline(-150,0)(-150,80)

\uput[90](-150,80){\small $F_2$}

\psline(-50,0)(-50,80)

\uput[90](-50,80){\small $F_{b-1}$}

\psline(0,0)(0,80)

\uput[90](0,80){\small $F_b$}

\pscircle*(-115,92){2}
\pscircle*(-100,92){2}
\pscircle*(-85,92){2}

\uput[270](-100,0){\small $\Gamma$}

\parametricplot{-30}{30}{-100 100 t 9 mul sin mul add 40 t add}

\uput[90](-100,40){\small $C$}

\psline{->}(25,50)(10,50)

\rput[lt]{0}(30,60){
\parbox{1.8in}{\small A tangency between $C$ and $F$ decreases $i_R(C,
F)$ by one but increases the total ramification index
of $p_C$ by one at the same time.}
}

\end{pspicture}

\caption{Application of \propref{PROP003}}
\label{FIG003}
\end{figure}

Second, it is also easy to verify that the proposition holds
when $n = k$ (we set $M_i = 0$ when $n = k$). When $n = k$, $R\isom
\P^k$, $F$ is a union of $b$ hyperplanes in $\P^k$ in very general
position and $L$ is a very general hypersurface of degree $a$ in
$\P^k$. Since $k < 3$ and we have the inductive hypothesis
that \thmref{THM003} holds in $\dim X < 3$,
\begin{equation}\label{E412}
2 g(C) - 2 + i_R(C, F\cup L) \ge (a+b - 2k) \deg C
\end{equation}
and \eqref{E222} follows.

So the proposition holds for $a = n - k$ or $n = k$.
We will also prove the proposition for $k = 1$.
These are the starting points of our induction. 
Next, we will use a degeneration
argument to bring down the value of $a$, $k$ or $n$ (the proof for the
case $k = 1$ is included in the following argument).

Let $\hat L$ be a very general member of the linear series
$|M+(a-1)\f|$ and $G$ be a very general member of the linear series
$|\f|$. By degenerating $L$ to $\hat L\cup G$, we will lower the value
of $a$, $k$ or $n$. Repeating this process, we will eventually reduce
the proposition to one of the cases we have already verified.

Let $Z = R\times\Delta$, $W_M = g^* M$, $W_F = g^* F$ and
$W_L\subset Z$ be a pencil in $|M+a\f|$ whose general fibers are
general members of the linear series and whose central fiber is $\hat
L\cup G$, where $g$ is the projection $Z\to R$.

Let $Y$ be a reduced irreducible flat family of curves over
$\Delta$ with the commutative diagram \eqref{E204}.
We assume that $\pi_* Y_t\cdot \f \ne 0$. Our goal is to prove that
\begin{equation}\label{E228}
2 g(Y_t) - 2 + i_Z(\pi_* Y_t, W) \ge (a+b - 2n) (\pi_* Y_t\cdot \f)
\end{equation}
for $t\ne 0$, where $W = W_M + W_F + W_L$.

If $k > 1$, we blow up $Z$ along $G$.
Let $f: \wt{Z}\to Z$ be the blowup map and $\wt{W}$
be the proper transform of $W$. As before, we
have the commutative diagram \eqref{E210} and the identity
\eqref{E211}. Let $E\subset \wt{Z}$ be the exceptional divisor of
$f$. Then the central fiber of $\wt{Z}$ is the union $R\cup E$ and
$R\cap E = G$.

Note that we do not blow up $Z$ if $k = 1$. If $k = 1$, we let $\wt{Z}
= Z$, $f$ be the identity map and $\wt{W} = W$ in the following
argument.

After resolving the indeterminacies of the rational map $Y\to \wt{Z}$,
we obtain the commutative diagram \eqref{E210}.

Again, we may apply \propref{PROP002} to $(\wt{Y}, \wt{Z}, \wt{W})$ and
obtain
\begin{equation}\label{E229}
\begin{split}
&\quad \quad \quad 2 g(Y_t) - 2 + i_Z(\pi_* Y_t, W)\\
&\ge \sum_{\genfrac{}{}{0pt}{}{\Gamma\subset E}{f_* \Gamma\ne 0}}
\mu_\Gamma \Phi_{\wt{Z},\wt{W}}(\Gamma)
+\sum_{\Gamma\not\subset E} \mu_\Gamma\Phi_{Z,W}(f_* \Gamma)
\end{split}
\end{equation}
where we sum over the irreducible components
$\Gamma\subset \wt{\pi}_* \wt{Y}_0$ satisfying that $f_* \Gamma \ne 0$
and $\Gamma\not\subset U_\wt{Y}$.
Here we need to exclude $\Gamma\subset U_\wt{Y}$ in the summations,
for reasons that will be clear later. 

Therefore, we just have to show
\begin{equation}\label{E231}
\Phi_{\wt{Z},\wt{W}}(\Gamma) \ge (a+b - 2n)(f_* \Gamma \cdot \f)
\end{equation}
for all $\Gamma\subset E$ and $f_*\Gamma \ne 0$ and
\begin{equation}\label{E230}
\Phi_{Z,W}(f_* \Gamma) \ge (a+b - 2n)(\Gamma \cdot \f)
\end{equation}
for all $\Gamma\subset R$ and $\Gamma\not\subset U_\wt{Y} \cup G$;
\eqref{E228} will follow easily.

Suppose that $\Gamma\subset E$ and $f_*\Gamma \ne 0$. Let us first
figure out what $E$ and $\wt{W}\cap E$ look like.

Obviously, $E = \P_{G} \E$ while $G\isom \P_{\P^{k-1}} \E^{(n-k)}$
is the projective tower of $n-k$ $\E$'s
over $\P^{k-1}$. Therefore, $E$ is the projective tower of
$n-k+1$ $\E$'s over $\P^{k-1}$, i.e., $E\isom \P_{\P^{k-1}}
\E^{(n-k+1)}$. 

It is not hard to see that
\begin{equation}\label{E402}
\wt{W}\cap E = \left(\bigcup_{j=1}^{n-k} M_j' \right)\cup
\left(\bigcup_{\beta=1}^b F_\beta'\right)\cup L'
\end{equation}
where
\begin{enumerate}
\item $M_j'\cap G = M_j\cap G$ for $1\le j\le n-k$ and $M_1',
M_2', ..., M_{n-k}'$ and $M_{n-k+1}' = G$ are the $n-k+1$
tautological divisors of $E \isom \P_{\P^{k-1}} \E^{(n-k+1)}$;
\item $F_\beta'\cap G = F_\beta\cap G$ and $F_\beta'$
is a very general member of
the linear series $\P H^0(\CO_E(\f'))$ for $\beta = 1, 2, ..., b$
with $\f'$ the pullback of the hyperplane divisor under the
projection $E\to \P^{k-1}$;
\item $L'\cap G = \hat L\cap G$ and $L'$ is a very general member of the
linear series $\P H^0(\CO_E(M' + a\f'))$ with $M' = M_1' + M_2' +
... + M_{n-k+1}'$.
\end{enumerate}
Thus
\begin{equation}\label{E236}
\begin{split}
\Phi_{\wt{Z},\wt{W}}(\Gamma) &= 2 g(\Gamma) - 2 +
i_E(\Gamma, (\wt{W}\cap E)\cup G)\\
&= 2 g(\Gamma) - 2 + i_E(\Gamma,  M'\cup F'\cup L')\\
&\ge (a + b - 2n) (\Gamma \cdot \f') = (a+b - 2n) (f_*\Gamma \cdot \f)
\end{split}
\end{equation}
by the inductive hypothesis, where $F' = F_1'+ F_2' + ... + F_b'$.
And \eqref{E231} follows.

Suppose that $\Gamma\subset R$, $\Gamma\not\subset U_\wt{Y}\cup G$ and
$p_* \Gamma = 0$, i.e.,
$\Gamma$ is contained in a fiber of $p: R\to \P^k$. Note that
$\Phi_{Z,W}(f_* \Gamma)$ could be negative in this case if we do not
exclude $\Gamma\subset U_\wt{Y}$. However, since $\Gamma\not\subset
U_\wt{Y}$,
\begin{equation}\label{E403}
\Phi_{Z,W}(f_* \Gamma)\ge 0
\end{equation}
by the definition of $U_\wt{Y}$. This is the reason that we need to
exclude those $\Gamma\subset U_\wt{Y}$.

Suppose that $\Gamma\subset R$, $p_*\Gamma \ne 0$
and $\Gamma\not\subset \hat L$. Then
\begin{equation}\label{E232}
\begin{split}
\Phi_{Z,W}(f_* \Gamma) &= 2 g(\Gamma) - 2 + i_R(\Gamma, W_0)\\
&= 2 g(\Gamma) - 2 + i_R(\Gamma, M \cup F \cup G\cup \hat L)\\
&\ge \left((a-1) + (b+1) - 2n\right) (\Gamma \cdot \f)\\
&= (a + b - 2n) (\Gamma \cdot \f)
\end{split}
\end{equation}
by the inductive hypothesis.

Suppose that $\Gamma\subset R$, $p_*\Gamma \ne 0$
and $\Gamma\subset \hat L$. Then
\begin{equation}\label{E233}
\begin{split}
\Phi_{Z,W}(f_* \Gamma)
&= 2 g(\Gamma) - 2 + i_R(\Gamma, M\cup F)\\
&= 2 g(\Gamma) - 2 + i_{\hat L}(\Gamma, (M\cup F)\cap \hat L).
\end{split}
\end{equation}
Let $\psi : X \to M_{n-k} = N$ be the projection from $X$ to $M_{n-k}$
where $X$ is regarded as the projectivization of $\E$ over $N = M_{n-k}$.
It is not hard to see that $N$ is the projective
tower of $n-k-1$ $\E$'s over $\P^k$, i.e.,
\begin{equation}\label{E238}
N \isom \P_{\P^k} \E^{(n-k-1)}.
\end{equation}
Let $F_\beta' = F_\beta\cap N$ for $\beta = 1,2,...,b$,
$M_j' = M_j \cap N$ for
$j=1,2,...,n-k-1$ and $L' = \hat L\cap N$. It is not hard to
see that
\begin{enumerate}
\item $F_\beta'$ is a very general member of the linear series
$\P H^0(\CO_N(\f'))$ for
$\beta = 1,2,...,b$, where $\f' = \f\cap N$ is the pullback of the
hyperplane divisor under the projection $N\to \P^{k}$;
\item $M_1', M_2', ..., M_{n-k-1}'$ are the $n-k - 1$
tautological divisors of $N \isom \P_{\P^k} \E^{(n-k -1)}$;
\item $L'$ is a very general member of the linear series
\begin{equation}\label{E237}
\P H^0\left(\CO_N(M'+ (a-2)\f'\right),
\end{equation}
where $M' = M_1' + M_2' + ... + M_{n-k - 1}'$.
\end{enumerate}
Since $\psi$ induces an isomorphism between $\hat L$ and $N$, we have
\begin{equation}\label{E234}
i_{\hat L}(\Gamma, (M\cup F)\cap \hat L)
= i_N(\psi_* \Gamma, M' \cup F' \cup L')
\end{equation}
where $F' = F_1' + F_2' + ... + F_b'$.
Apply the inductive hypothesis to
$(N, M' \cup F' \cup L')$ and we obtain
\begin{equation}\label{E235}
\begin{split}
\Phi_{Z,W}(f_* \Gamma)
& = 2 g(\Gamma) - 2 + i_N(\psi_* \Gamma, M' \cup F' \cup L')\\
& \ge \left((a-2) + b - 2(n-1)\right) (\psi_* \Gamma \cdot \f')\\
&= (a+b - 2n) (\Gamma \cdot F).
\end{split}
\end{equation}
Combine \eqref{E403}, \eqref{E232} and \eqref{E235} and we obtain
\eqref{E230}.
\end{proof}

\section{Proof of \thmref{THM003}}\label{SEC004}

We will carry out the proof of \thmref{THM003} in two
steps as outlined in \ssecref{SEC003001}: we will first degenerate
$D$ and then $S$.

One of our main tools is \propref{PROP002}, which
will be proved at the end of this section.

\subsection{Degeneration of $D$}\label{SEC004001}

Without the loss of generality, we assume that $D_i\ne 0$ for all
$i\in I$ and $I_\gamma \ne \emptyset$ for all $\gamma\in \CA\cup\CB$.
Hence $\dim S \le \dim X - |\CA|$. Our inductive hypothesis is:

\begin{hyp}\label{HYP001}
The theorem holds in $\dim S < n$.
\end{hyp}

Of course, to start the induction, we have to verify the following.

\begin{claim}\label{CLM001}
The theorem holds in $\dim S = 1$.
\end{claim}

\begin{proof}
By Bertini, $S$ meets $D$ transversely. Therefore, $i_X(S, D) = D\cdot
S$ and
\begin{equation}\label{E001}
2 g(S) - 2 + i_X(S, D) = \left(K_X + \sum_{i\in I} D_i\right) S.
\end{equation}
Let $C$ be an irreducible component of $S$. Since $S$ is cut by
general members of BPF linear systems, the monodromy actions on the
irreducible components of $S$ are transitive, which implies that any
two irreducible component of $S$ are numerically
equivalent. Therefore,
\begin{equation}\label{E004}
2 g(C) - 2 + i_X(C, D) = \left(K_X + \sum_{i\in I} D_i\right) C.
\end{equation}
Since $\{ D_i \}_{i\in I}$ satisfies the condition $\C_X(\varphi, \dim
X - 2, 0)$, the curves $C$ satisfying
\begin{equation}\label{E005}
\left(K_X + \sum_{i\in I} D_i\right) C < \varphi(C)
\end{equation}
cannot cover $X$ by \lemref{LEM001}. However, if $C$ is an irreducible
component of $S$, the curves in the numerical class of $C$ obviously
covers $X$. Consequently, we necessarily have
\begin{equation}\label{E006}
\left(K_X + \sum_{i\in I} D_i\right) C \ge \varphi(C)
\end{equation}
and \eqref{E308} follows.
\end{proof}

If $|I_\beta| = 1$ for every $\beta\in \CB$, then we proceed directly to
\ssecref{SEC004002}.
Otherwise, let $|I_\gamma| > 1$ for some $\gamma\in \CB$.
We will degenerate $D_{I_\gamma}$ to $\sum_{i\in I_\gamma} D_i$. 
However, as in the special case carried out in
last section, instead of degenerating $D_{I_\gamma}$
directly to $\sum_{i\in I_\gamma} D_i$,
we will degenerate $D_{I_\gamma}$ one component at a time. 

Pick an arbitrary $\kappa\in I_\gamma$ and let $\hat I_{\gamma} =
I_\gamma \backslash \{\kappa\}$ and $D_{\hat I_\gamma}$ be a
very general member of the linear series $\P\L_{\hat I_{\gamma}}$. We
will degenerate $D_{I_\gamma}$ to the union $D_{\hat I_{\gamma}} \cup
D_\kappa$ and keep doing this until $|I_\beta| = 1$ for every
$\beta\in \CB$.

The argument for each step of degeneration, i.e., the
degeneration of $D_{I_\gamma}$ to $D_{\hat I_\gamma}\cup D_\kappa$, has been
amply illustrated by the special case we did in the previous
section. Basically, it involves the construction of a fan over $G =
D_\kappa$ and then the proof of a statement similar to \propref{PROP004}.

Let $Z = X\times\Delta$ and $\D_{I_\gamma}\subset Z$
be a pencil in the linear
series $\P\L_{I_\gamma}$ whose general fibers are general members of
the linear series and whose central fiber is
$D_{\hat I_{\gamma}} \cup D_\kappa$.

Let $\D_{I_\tau} = p^* D_{I_\tau}$ for $\tau\in \CA\cup \CB$ and
$\tau\ne \gamma$, $\S = \cap_{\alpha\in\CA} \D_{I_\alpha}$ and $\D =
\sum_{\beta\in \CB} \D_{I_\beta}$,
where $p: Z\to X$ is the projection from $Z$ to $X$.

Let $Y$ be a reduced flat family of curves with the commutative
diagram \eqref{E204}. We assume that $\pi(Y)\subset \S$.
Our goal is to prove
\begin{equation}\label{E018}
2 g(Y_t) - 2 + i_Z(\pi_* Y_t, \D) \ge \varphi(\pi_* Y_t)
\end{equation}
for $t\ne 0$ under the inductive hypothesis:
\begin{equation}\label{E002}
2 g(C) - 2 + i_X (C, \D_0) \ge \varphi(C)
\end{equation}
for all reduced curves $C\subset S$.

Let $f: \wt{Z}\to Z$ be the blowup of $Z$ along $G = D_\kappa$ and let
$R$ be the exceptional divisor. Then $\wt{Z}_0 = X\cup R$ and $X\cap
R = G$. Let $\wt{\D}$ and $\wt{\S}$ be the proper transforms of $\D$
and $\S$ under $f$, respectively.

After resolving the indeterminacies of the rational map $Y\to \wt{Z}$,
we obtain the commutative diagram \eqref{E210}.

By \propref{PROP002}, we have
\begin{equation}\label{E019}
\begin{split}
&\quad\quad\quad 2 g(Y_t) - 2 + i_Z(\pi_* Y_t, \D)\\
&\ge \sum_{\genfrac{}{}{0pt}{}{\Gamma\subset R}{f_* \Gamma\ne 0}}
\mu_\Gamma \Phi_{\wt{Z},\wt{\D}}(\Gamma) + 
\sum_{\Gamma\not\subset R} \mu_\Gamma \Phi_{Z,\D}(f_* \Gamma)
\end{split}
\end{equation}
where we sum over the irreducible components
$\Gamma\subset \wt{\pi}_* \wt{Y}_0$ not contracted by $f$. 

Let $P = D_{\hat I_{\gamma}}$. By the definition of $\Phi$, we have
\begin{align}\label{E413}
\Phi_{\wt{Z},\wt{\D}}(\Gamma) &= 2 g(\Gamma) - 2 +
i_R(\Gamma, G \cup (\wt{\D}\cap R)) &&\text{for } \Gamma\subset R,\\
\label{E414}
\Phi_{Z,\D}(f_* \Gamma) &= 2 g(\Gamma) - 2 + i_Z(f_* \Gamma, \D)
&&\text{for } \Gamma\subset P, \text{ and}\\
\label{E415}
\Phi_{Z,\D}(f_* \Gamma) &= 2 g(\Gamma) - 2 + i_X (f_* \Gamma, \D_0)
&&\text{for } \Gamma\not\subset R\cup P.
\end{align}
Thus it is enough to show the following:
\begin{equation}\label{E020}
2 g(\Gamma) - 2 + i_R(\Gamma, G \cup (\wt{\D}\cap R))\ge \varphi(f_* \Gamma)
\end{equation}
for any $\Gamma\subset \wt{\S}\cap R$ and $f_* \Gamma\ne 0$,
\begin{equation}\label{E021}
2 g(\Gamma) - 2 + i_Z(\Gamma, \D) \ge \varphi(\Gamma)
\end{equation}
for any $\Gamma\subset S\cap P$ and
\begin{equation}\label{E022}
2 g(\Gamma) - 2 + i_X (\Gamma, \D_0) \ge
\varphi(\Gamma)
\end{equation}
for any $\Gamma\subset S$ and $\Gamma\not\subset P\cup G$.

It is easy to see that \eqref{E022} follows directly from the inductive
hypothesis \eqref{E002}. And \eqref{E021} follows from \Hypref{HYP001}
by the following argument. Let $\hat D = D - D_{I_\gamma}$,
$\hat S = S\cap P = \cap_{\alpha\in \CA} D_{I_\alpha}\cap D_{\hat
I_\gamma}$. By \lemref{LEM001}, $\{ D_i: i\in I, i\ne \kappa \}$
satisfies $\C_X(\varphi, \dim X - 2, \dim \hat S - 1)$.
Therefore,
\begin{equation}\label{E003}
2 g(C) - 2 + i_X(C, \hat D) \ge \varphi(C)
\end{equation}
for all reduced curves $C\subset \hat S$ by \Hypref{HYP001}, since
$\dim \hat S \le \dim S - 1$. Then \eqref{E021} follows by observing
$i_Z(\Gamma, \D) = i_X(\Gamma, \hat D)$ for $\Gamma\subset\hat
S$.

It remains to justify \eqref{E020}. This leads to the following
proposition similar to \ref{PROP004}, from which \eqref{E020} follows.

\begin{prop}\label{PROP001}
Let $X$ be a smooth projective variety, $\{ \P \L_i\}_{i\in I}$ be a
finite set of BPF linear systems on
$X$, $D_i$ be a very general member of $\P\L_i$ and
$\E_i = \CO\oplus \CO(D_i)$ for $i\in I$.

Let $J\subset I$ and $R = \P_X \E^{(J)} = \P_X (\E_j)_{j\in J}$
be the projective tower of
$\E_j$ for $j\in J$, i.e., the fiber product 
\begin{equation}\label{E401}
\prod_{j\in J} \P_X \E_j
\end{equation}
over $X$. Let $M_j$ be the pullback of the tautological
divisor of $\P_X \E_j$ under the projection $p_j: R\to \P_X \E_j$ and $M =
\sum_{j\in J} M_j$.

Let
\begin{equation}\label{E026}
I = \left(\bigsqcup_{\alpha\in\CA} I_\alpha\right) \sqcup
\left(\bigsqcup_{\beta\in \CB} I_\beta\right) \sqcup
I_\gamma\sqcup J
\end{equation}
be a partition of $I$. Let $D_{I_\tau}$ be a
very general member of $\P \L_{I_\tau}$ and $F_\tau = p^*
D_{I_\tau}$ for $\tau\in\CA\cup \CB$,
where $p: R\to X$ is the natural projection.
Let $L$ be a very general member of the linear series
\begin{equation}\label{E027}
\P(\tensor_{j\in J} (\BC\oplus \L_j)) \otimes \L_{I_\gamma}
\subset \P H^0(\CO_R(M + \sum_{j\in J\cup I_\gamma} p^* D_j)),
\end{equation}
$S = \cap_{\alpha\in \CA} D_{I_\alpha}$, $\S = \cap_{\alpha\in \CA}
F_\alpha$, $D = \sum_{\beta\in\CB} D_{I_\beta}$ and
$F = \sum_{\beta\in\CB} F_\beta$. Here we use the identification between
the linear systems on $X$ and those on $R$, as given in \eqref{E301}.

Let $\varphi: N_1(X)\to \BR$ be a subadditive function on $N_1(X)$
such that $\{ D_i \}_{i\in I\backslash J}$ satisfies the condition
$\C_X(\varphi, \dim X - 2, \dim \S - 1)$ on $X$. Suppose that
$0 < \dim S < n$. Then
\begin{equation}\label{E007}
2 g(C) - 2 + i_R(C, M\cup F\cup L) \ge \varphi(p_* C)
\end{equation}
for all reduced irreducible curves $C\subset \S$ satisfying $p_* C \ne
0$.
\end{prop}

\begin{proof}
First, let us verify the case that $I_\gamma = \emptyset$.
In this case, since $\dim S < n$,
\begin{equation}\label{E008}
2 g(\Gamma) - 2 + i_X(\Gamma, D) \ge \varphi(\Gamma)
\end{equation}
for all reduced curves $\Gamma\subset S$ by \Hypref{HYP001}.
Let $\Gamma$ be the reduced image of $C$ under the projection $p$ and
let $p_C: C\to \Gamma$ be the restriction of $p$ to $C$.
As in the proof of \propref{PROP004}, we may apply \propref{PROP003} to
the map $p_C$ to obtain (see \figref{FIG003}):
\begin{equation}\label{E010}
2 g(C) - 2 + i_R(C, F) \ge (\deg p_C)(2 g(\Gamma) - 2 + i_X(\Gamma,
D)).
\end{equation}
Then \eqref{E007} follows from \eqref{E008} and \eqref{E010}.

Second, the proposition obviously holds when $J = \emptyset$, i.e., $X = R$
by \Hypref{HYP001}.

So far we have verified the proposition for
$I_\gamma = \emptyset$ or $X = R$. We will also prove the proposition
for $\dim S = 1$.
Next we will try to reduce it to one of these cases by degenerating
$L$ (the proof for $\dim S = 1$ is included in the following argument).

Pick an arbitrary $\kappa\in I_\gamma$ and let $I_\gamma = \hat
I_\gamma\sqcup \{\kappa\}$ and $\hat L$ be a very general member of
the linear series
\begin{equation}\label{E012}
\P(\tensor_{j\in J} (\BC\oplus \L_j)) \otimes \L_{\hat I_\gamma}
\subset \P H^0(\CO_R(M + \sum_{j\in J\cup {\hat I_\gamma}} p^* D_j)).
\end{equation}
By degenerating $L$ to $\hat L \cup D_\kappa$, we will lower the
value of $|I_\gamma|$, $\dim X$ or $\dim R - \dim X$. So repeating
this process, we will eventually reduce the proposition to one of the
cases we have proved.

Let $Z = R\times\Delta$, $W_\beta = g^* F_\beta$ for $\beta\in
\CB$ and $W_\gamma\subset Z$ be a pencil in the linear series
\eqref{E027} whose general fibers are general members of the
linear series and whose central fiber is $\hat L \cup D_\kappa$, where
$g$ is the projection $Z\to R$. Let $W = W_\gamma + \sum_{\beta\in\CB}
W_\beta + \sum_{j\in J} g^* M_j$.

Let $Y$ be a reduced irreducible flat family of curves with the commutative
diagram \eqref{E204}. We assume that $\pi(Y)\subset \S$ and
$p_* \pi(Y) \ne 0$. Our goal is to prove
\begin{equation}\label{E013}
2 g(Y_t) - 2 + i_Z(\pi_* Y_t, W) \ge \varphi(\pi_* Y_t).
\end{equation}

If $\dim S > 1$, we blow up $Z$ along $G = D_\kappa$.
Let $f: \wt{Z}\to Z$ be the blowup map and let
$E$ be the exceptional divisor. Then $\wt{Z}_0 = R\cup E$ and $R\cap
E = G$. Let $\wt{W}_\beta$, $\wt{W}_\gamma$ and $\wt{W}$ be
the proper transforms of $W_\beta$, $W_\gamma$ and $W$, respectively.

Note that we do not blow up $Z$ if $\dim S = 1$. If $\dim S = 1$, we
let $\wt{Z} = Z$, $f$ be the identity map and $\wt{W} = W$ in the
following argument.

After resolving the indeterminacies of the rational map $Y\to \wt{Z}$,
we obtain the commutative diagram \eqref{E210}.

By \propref{PROP002},
\begin{equation}\label{E014}
\begin{split}
&\quad\quad\quad 2 g(Y_t) - 2 + i_Z(\pi_* Y_t, W)\\
&\ge \sum_{\genfrac{}{}{0pt}{}{\Gamma\subset E}{f_* \Gamma\ne 0}}
\mu_\Gamma \Phi_{\wt{Z}, \wt{W}}(\Gamma)
+ \sum_{\Gamma\not\subset E} \mu_\Gamma \Phi_{Z, W}(f_* \Gamma)
\end{split}
\end{equation}
where we sum over the irreducible components $\Gamma\subset \wt{\pi}_*
\wt{Y}_0$ satisfying that $f_* \Gamma\ne 0$ and $\Gamma\not\subset
U_\wt{Y}$.
Here we need to exclude $\Gamma\subset U_\wt{Y}$ due to the presence
of the components $\Gamma\subset R$ with $p_* \Gamma = 0$.

By the definition of $\Phi$, we have
\begin{align}\label{E011}
\Phi_{\wt{Z},\wt{W}}(\Gamma) &=
2 g(\Gamma) - 2 + i_E(\Gamma, G\cup (\wt{W}\cap E)) &&\text{for }
\Gamma\subset E,\\
\label{E500}
\Phi_{Z,W}(f_* \Gamma) &= 2 g(\Gamma) - 2 + i_Z(f_* \Gamma, W) &&
\text{for } \Gamma\subset \hat L, \text{ and}\\
\label{E501}
\Phi_{Z,W}(f_* \Gamma) &= 2 g(\Gamma) - 2 + i_R(\Gamma, W_0)
&& \text{for } \Gamma\not\subset E\cup \hat L.
\end{align}

Observe that $E$ is actually the projective tower $\P_G \E^{(J)}
\times_G \P_G \E_\kappa$ over $G$ and 
\begin{equation}\label{E015}
\wt{W}\cap E = \left(\bigcup_{j\in J} M_j'\right) \cup
\left(\bigcup_{\beta\in\CB} F_\beta'\right) \cup L'
\end{equation}
where
\begin{enumerate}
\item $M_j'\cap G = M_j\cap G$ and $M_j'$ and $G$ are the
tautological divisors of $E$ over $G$ for $j\in J$;
\item $F_\beta'\cap G = F_\beta\cap G$ and $F_\beta' = f^*(D_{I_\beta}
\cap G)$;
\item $L'\cap G = \hat L\cap G$ and $L'$ is a very general member of
the linear series
\begin{equation}\label{E016}
\P(\tensor_{j\in J} (\BC\oplus \L_j)) \otimes (\BC \oplus \L_\kappa)
\otimes \L_{\hat I_\gamma}
\end{equation}
with $\L_i$ restricted to $G$.
\end{enumerate}
Note that we have $\dim G < \dim X$ and $\{ D_i: i\in I, i\ne \kappa
\}$ satisfies the condition $\C_G(\varphi, \dim G - 2, \dim \S - 1)$
when restricted to $G$.
Therefore, we may apply the inductive hypothesis to $(E, M'\cup F' \cup
L')$ and obtain
\begin{equation}\label{E017}
\Phi_{\wt{Z}, \wt{W}}(\Gamma) = 2 g(\Gamma) - 2 +
i_E(\Gamma, M'\cup F' \cup L') \ge \varphi(p_* f_* \Gamma)
\end{equation}
for any $\Gamma\subset E$, $f_*\Gamma \ne 0$ and $f(\Gamma)\subset S$,
where $M' = \sum_{j\in J} M_j' + G$ and $F' = \sum_{\beta\in\CB}
F_\beta'$.

Suppose that $\Gamma\subset R$, $\Gamma\not\subset U_\wt{Y}\cup G$ and
$p_* \Gamma = 0$,
i.e., $\Gamma$ is contained in a fiber of $p$.
Note that $\Phi_{Z, W}(f_* \Gamma)$ could be negative if we do not
exclude $\Gamma \subset U_\wt{Y}$.
But since $\Gamma \not\subset U_\wt{Y}$, it is easy to check that
\begin{equation}\label{E023}
\Phi_{Z, W}(f_* \Gamma) \ge 0
\end{equation}
by the definition of $U_\wt{Y}$.
This is the reason that we need to exclude those $\Gamma \subset U_\wt{Y}$.

Suppose that $\Gamma\subset R$, $p_* \Gamma\ne 0$ and
$\Gamma\not\subset \hat L$. Then it follows directly from the
inductive hypothesis that
\begin{equation}\label{E024}
\Phi_{Z, W}(f_* \Gamma) = 2 g(\Gamma) - 2 + i_R(\Gamma, W_0)
\ge \varphi(p_* \Gamma)
\end{equation}
for any $\Gamma\subset \S$, $p_* \Gamma\ne 0$ and
$\Gamma\not\subset \hat L$.

Suppose that $\Gamma\subset R$, $p_* \Gamma\ne 0$ and
$\Gamma\subset \hat L$. We do exactly the same thing as we did in
the proof of \propref{PROP004} by projecting $R$ to $M_\tau$ for some
$\tau \in J$. The projection $\psi: R\to M_\tau = N$ induces an isomorphism 
between $\hat L$ and $N$. It is not hard to see that $N$ is the projective
tower $\P_X\E^{(\hat J)} = \P_X(\E_j)_{j\in \hat J}$
over $X$ with $\hat J = J\backslash \{\tau\}$.

Let $F_\beta' = F_\beta\cap N$ for $\beta\in\CB$, $M_j' = M_j \cap N$ for
$j\in \hat J$ and $L' = \hat L\cap N$. It is not hard to
see that
\begin{enumerate}
\item $F_\beta' = p_N^* D_{I_\beta}$, where $p_N$ is the projection $N\to X$;
\item $M_j'$ are the tautological divisors of $N$ for $j\in\hat J$;
\item $L'$ is a very general member of the linear series
\begin{equation}\label{E025}
\P(\tensor_{j\in \hat J} (\BC\oplus \L_j)) \otimes \L_{\hat I_\gamma}
\end{equation}
on $N$.
\end{enumerate}

Since $\psi$ induces an isomorphism between $\hat L$ and $N$, we have
\begin{equation}\label{E028}
i_{\hat L}(\Gamma, (M \cup F)\cap \hat L)
= i_N(\psi_* \Gamma, M' \cup F' \cup L')
\end{equation}
where $M' = \sum_{j\in\hat J} M_j'$ and $F' = \sum_{\beta\in\CB}
F_\beta'$. Let $\S' = \S\cap \hat L$.
Since $\dim N - \dim X < \dim R - \dim X$ and $\{ D_i: i\in I, i\ne \kappa
\}$ satisfies the condition $\C_X(\varphi, \dim X - 2, \dim \S' - 1)$,
we may apply the inductive hypothesis to
$(N, M' \cup F' \cup L')$ and obtain
\begin{equation}\label{E029}
2 g(\psi_* \Gamma) - 2 + i_N(\psi_* \Gamma, M' \cup F' \cup L')
\ge \varphi(p_* \Gamma)
\end{equation}
for all $\Gamma\subset \S'$ and $p_*\Gamma \ne 0$.
Therefore,
\begin{equation}\label{E030}
\begin{split}
\Phi_{X,W}(f_*\Gamma)
&= 2 g(\Gamma) - 2 + i_{\hat L}(\Gamma, (M \cup F)\cap \hat L)\\
&= 2 g(\psi_* \Gamma) - 2 + i_N(\psi_* \Gamma, M' \cup F' \cup L')
\ge \varphi(p_* \Gamma)
\end{split}
\end{equation}
for all $\Gamma\subset \S'$ and $p_*\Gamma \ne 0$.
\end{proof}

Combining \eqref{E014},
\eqref{E017}, \eqref{E023}, \eqref{E024} and \eqref{E030}, we arrive
at \eqref{E013}.

\subsection{Degeneration of $S$}\label{SEC004002}

After we have reduced the theorem to the case that $|I_\beta| = 1$ for all
$\beta\in\CB$ by degenerating $D$, we will finish the proof by
degenerating $S$.

For every $\alpha\in \CA$, let $Z_\alpha\subset X\times\Delta$ be a
pencil in the linear series $\P\L_{I_\alpha}$ whose general fibers
are general members of the linear series and whose central fiber is
the union $\cup_{i\in I_\alpha} D_i$. Let $Z = \cap_{\alpha\in\CA}
Z_\alpha$ and $W = p^* D$, where $p$ is the projection $Z\to X$.

Let $Y$ be a reduced flat family of curves with the commutative
diagram \eqref{E204}. Our goal is to prove that
\begin{equation}\label{E031}
2 g(Y_t) - 2 + i_Z(\pi_* Y_t, W) \ge \varphi(\pi_* Y_t).
\end{equation}

The central fiber of $Z$ is a union of projective varieties of simple
normal crossing with components $\cap_{\alpha\in\CA} D_{i_\alpha}$,
where $i_\alpha\in I_\alpha$. By \propref{PROP002},
\begin{equation}\label{E032}
2 g(Y_t) - 2 + i_Z(\pi_* Y_t, W) \ge \sum_{\Gamma\subset\pi_* Y_0}
\mu_\Gamma \Phi_{Z,W}(\Gamma).
\end{equation}
Obviously, in order to show \eqref{E031}, it suffices to show
\begin{equation}\label{E033}
\Phi_{Z,W}(\Gamma) \ge \varphi(\Gamma)
\end{equation}
for every irreducible component $\Gamma\subset \pi_* Y_0$.

Let $J = \{ j\in I: \Gamma\subset D_j \}$,
$M_\alpha = I_\alpha\backslash J$
and
$D_J = \cap_{j\in J} D_j$. Then
\begin{equation}\label{E034}
\Phi_{Z,W}(\Gamma) \ge 2 g(\Gamma) - 2 + i_{D_J}(\Gamma, (\cup_{i\in
I\backslash J} D_i)\cap D_J, Q\cap \partial D_J)
\end{equation}
where $Q$ is the singular locus of $Z$ and
\begin{equation}\label{E035}
\partial D_J = D_J\cap \left(\bigcup_{\alpha\in\CA} \bigcup_{i\in
M_\alpha} D_i\right).
\end{equation}
Let us assume that $\dim D_J\ge 2$ since \eqref{E033} follows
from \lemref{LEM001} if $\dim D_J = 1$ by the argument for
the claim \ref{CLM001}.

The singular locus $Q\cap \partial D_J$ can be described as follows:
let $F_\alpha$ be a general member of the pencil $Z_\alpha$ and then
\begin{equation}\label{E036}
Q\cap \partial D_J = D_J\cap \left(\bigcup_{\alpha\in\CA}
\bigcup_{i\in M_\alpha} (D_i\cap F_\alpha)\right).
\end{equation}

Since $\{ D_i \}_{i\in I}$ satisfies the condition $\C_X(\varphi, \dim
X - 2, \dim Z_0 - 1)$,
$\{ D_i \}_{i\in I\backslash J}$ satisfies the
condition $\C_{D_J}(\varphi, \dim D_J - 2, \dim D_J - 1)$ when
restricted to $D_J$ by \lemref{LEM001}.
Then it follows from \thmref{THM002} that
\begin{equation}\label{E037}
2 g(\Gamma) - 2 + i_{D_J}(\Gamma, (\cup_{i\in
I\backslash J} D_i)\cap D_J, Q\cap \partial D_J) \ge \varphi(\Gamma)
\end{equation}
for all reduced irreducible curves $\Gamma$ that are not contained in
$\cup_{\alpha\in\CA} F_\alpha$. 

Note that 
$F_\alpha$ is a general
member of the linear series $\P\L_{I_\alpha}$. Fix $Q$ and we consider
the linear subseries $\f_Q$ of $\P\L_{I_\alpha}$ given by
\begin{equation}\label{E038}
\f_Q = \{ F_\alpha\in \P\L_{I_\alpha}: F_\alpha\cap (\cup_{i\in
M_\alpha} D_i) = Q\cap (\cup_{i\in M_\alpha} D_i)\}.
\end{equation}
Since $M_\alpha\subsetneq I_\alpha$, the base locus of $\f_Q$ is
exactly $Q\cap (\cup_{i\in M_\alpha} D_i)$. Therefore, it is possible
to choose $F_\alpha\in \f_Q$ such that $\Gamma\not\subset
F_\alpha$. This is true for every $\alpha\in\CA$. So
it is possible to choose $\{F_\alpha\}_{\alpha\in\CA}$
with $Q$ fixed such that
$\Gamma\not\subset\cup_{\alpha\in\CA} F_\alpha$. Then \eqref{E037}
follows.

\subsection{Proof of \propref{PROP002}}\label{SEC004003}

The proof of \propref{PROP002} is almost the same as that of
\cite[Theorem 1.17]{C2}. First let us quote two lemmas in \cite{C2}
(Lemma 2.3 and 4.1).

\begin{lem}\label{LEM005}
Let \(X\subset \Delta_{x_1 x_2 ... x_r}^{r}\times \Delta\) be the
hypersurface given by $x_1x_2...x_n = t$ for some $n \le r$,
where \(\Delta_{x_1 x_2 ... x_r}^{r}\)
is the $r$-dimensional polydisk parameterized by
\((x_1, x_2, ..., x_r)\) and $\Delta$ is the disk parameterized by $t$. 
Let \(X_0 = D = \cup_{k=1}^n D_k\) with \(D_k = \{ x_k = t = 0 \}\).

Let
$Y$ be a flat family of curves over $\Delta$ with the commutative diagram
\eqref{E212}. Suppose that \(\pi_* Y_0\ne 0\).
Then for each \(D_k\) there exists a curve \(\Gamma' \subset
\pi(Y_0)\) with \(\Gamma' \subset D_k\).
\end{lem}

\begin{lem}\label{LEM006}
Let \(X\subset \Delta_{x_1 x_2 ... x_r}^{r}\times \Delta\) be the
hypersurface given by \(x_1x_2...x_n = t f(t, x_1, x_2, ..., x_r)\),
where \(n < r\), \(f(0, 0, 0, ..., 0) = 0\) and
\(f(t, x_1, x_2, ..., x_r) \not\equiv 0\) along
\(x_1 = x_2 = ... = x_n = t = 0\). Let \(X_0 = D = \cup_{k=1}^n D_k\)
where \(D_k = \{ x_k = t = 0 \}\). And let \(Q\) be the singular locus
of \(X\), i.e., \(Q\) is cut out on \(X\) by \(f(t, x_1, x_2, ...,
x_r) = 0\).

Let
$Y$ be a flat family of curves over $\Delta$ with the commutative diagram
\eqref{E212}. Suppose that there exists \(J\subset \{1,2,...,n\}\)
and a reduced irreducible curve \(\Gamma\subset \pi(Y_0)\)
such that \(\Gamma\subset
D_J\), \(\Gamma\not\subset\partial D_J\) and \(i_{D_J}(\Gamma,
\partial D_J, Q \cap \partial D_J) > 0\). Then there exists a
curve \(\Gamma'\subset \pi(Y_0)\) such that \(\Gamma' \subset
\cup_{k\not\in J} D_k\).
\end{lem}

Please see \cite{C2} for the proofs of these two lemmas.

\begin{defn}\label{DEF006}
Let $X, D$ and $P$ be given as in \defnref{DEF002} and let $f: C\to X$
be a proper map from a curve $C$ to $X$. We define
\begin{equation}\label{E040}
i_X(C, D, P)_f = i_X(f_* C, D, P).
\end{equation}
When there is no confusion on what $f$ is, we just write $i_X(C,D,P)_f$
as $i_X(C, D, P)$.
For a point $p\in C$, we define $i_{X, p}(C, D, P)$ by choosing an
open neighborhood (analytic or etale) $U$ of $f(p)$ and letting
\begin{equation}\label{E041}
i_{X, p}(C, D, P) = i_U(V, D\cap U, P\cap U)
\end{equation}
where $V$ is the connected component of $f^{-1}(U)$ that contains the
point $p$. Obviously, if $C$ is smooth,
\begin{equation}\label{E042}
i_X(C,D,P) = \sum_{p\in C} i_{X, p}(C, D, P).
\end{equation}
Finally, we define $I_X(C, D, P)\subset C$ to be the set
\begin{equation}\label{E043}
I_X(C, D, P) = \{ p\in C: i_{X, p}(C,D,P) \ne 0 \}.
\end{equation}
\end{defn}

Now we are ready to prove \propref{PROP002}.

Without the loss of generality, let us assume that $Y$ is irreducible
and birational over the the image $\pi(Y)$, which meets $W$ properly.
Furthermore, we may normalize $Y$ and apply 
semistable reduction to the map $Y\to X$. 
In the end, we may assume that $\pi: Y\to X$ is
a family of semistable maps with
marked points $Y_t\cap \pi^*(W)$ for $t\ne 0$. More
specifically, we assume that the following holds:
\begin{enumerate}
\item $Y$ is smooth and $Y_0$ is nodal;
\item $Y_t\cap \pi^*(W)$ extends to disjoint sections of the
fibration $Y\to \Delta$, i.e., the flat limit
$\lim_{t\to 0} (Y_t\cap \pi^*(W))$ consists of $i_X(\pi_* Y_t, W)$ distinct
points lying on the nonsingular locus of $Y_0$;
\item $\pi: Y\to X$ is minimal with respect to these properties.
\end{enumerate}
For each component $\Gamma\subset Y_0$, we define $\sigma(\Gamma)$
to be the set of points on $\Gamma$ that consists of the
marked points $\Gamma\cap \lim_{t\to 0}
(Y_t\cap \pi^*(W))$ and the intersections between $\Gamma$ and
the other components of $Y_0$. The basic
observation is (see \cite{C2}):
\begin{equation}\label{E044}
2 g(Y_t) - 2 + i_X(\pi_* Y_t, W) = \sum_{\Gamma\subset Y_0}
(2p_a(\Gamma) - 2 + |\sigma(\Gamma)|)
\end{equation}
where $p_a(\Gamma)$ is the arithmetic genus of $\Gamma$ and we
sum over all irreducible components $\Gamma\subset Y_0$. Therefore,
we will be done with \eqref{E213} as long as we can prove
\begin{equation}\label{E045}
2p_a(\Gamma) - 2 + |\sigma(\Gamma)| \ge (\deg \pi_\Gamma) \Phi_{X,W}(F)
\end{equation}
for each irreducible component $\Gamma\subset Y_0$, where
$F = \supp(\pi(\Gamma))$, $\pi_\Gamma: \Gamma\to F$ is the
restriction of $\pi$ to $\Gamma$ and we let $\Phi_{X,W}(F) = 0$ if $F
= \text{pt}$, i.e., $\Gamma$ is contractible under $\pi$.

Obviously, if $\Gamma$ is contractible under $\pi$,
\begin{equation}\label{E404}
2p_a(\Gamma) - 2 + |\sigma(\Gamma)| \ge 0
\end{equation}
due to the minimality of $\pi: Y\to X$.

Suppose that $\Gamma$ is noncontractible under $\pi$.
Let $\rho(\Gamma)\subset \Gamma$ be the set
\begin{equation}\label{E407}
\rho(\Gamma) =  I_X(\Gamma, W\cup D, Q\cap \partial D_J)
\end{equation}
where $J = \{j\in I: F\subset D_j\}$.
By \lemref{LEM005} and \ref{LEM006} and using the same argument as in
\cite{C2}, we can show that
\begin{equation}\label{E046}
\rho(\Gamma) \subset \sigma(\Gamma).
\end{equation}
And by \propref{PROP003},
\begin{equation}\label{E047}
\begin{split}
&\quad\quad 2 g(\Gamma) - 2 + |\rho(\Gamma)|\\
&\ge (\deg \pi_\Gamma)(2 g(F) - 2 + i_X(F, W\cup D, Q\cap \partial
D_J)).
\end{split}
\end{equation}
Then \eqref{E045} follows.

Let us consider the modified statement that excludes the components
$\Gamma\subset U_Y$ from the RHS of \eqref{E213}. Let $Z = \pi^{-1}(U_Y)$.

If $Z = Y_0$, there is nothing to prove since $U_Y =
\supp(\pi_* Y_0)$ in this case. Otherwise, we write
the RHS of \eqref{E044} as
\begin{equation}\label{E048}
\begin{split}
&\quad \quad \sum_{\Gamma\subset Y_0}
(2p_a(\Gamma) - 2 + |\sigma(\Gamma)|)\\
&= \sum_{\Gamma\not\subset Z} (2 p_a(\Gamma) - 2 + |\rho(\Gamma)|)\\
&+ \sum_{\Gamma\not\subset Z} (|\sigma(\Gamma)| - |\rho(\Gamma)|)
+ \sum_{\Gamma\subset Z} (2p_a(\Gamma) - 2 + |\sigma(\Gamma)|).
\end{split}
\end{equation}
We have already proved that
\begin{equation}\label{E405}
2 g(\Gamma) - 2 + |\rho(\Gamma)|\ge (\deg \pi_\Gamma) \Phi_{X,W}(F)
\end{equation}
for $\Gamma\not\subset Z$. So it suffices to show that
\begin{equation}\label{E406}
\sum_{\Gamma\not\subset Z} (|\sigma(\Gamma)| - |\rho(\Gamma)|)
+ \sum_{\Gamma\subset Z} (2p_a(\Gamma) - 2 + |\sigma(\Gamma)|) \ge 0.
\end{equation}

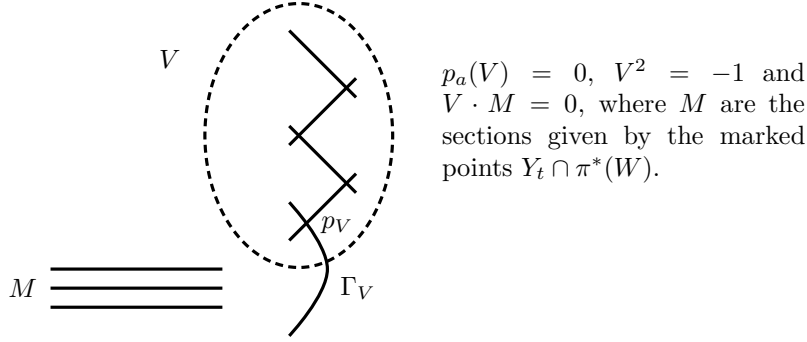
\begin{figure}[ht]
\centering
\begin{pspicture}(-230,-45)(230,140)

\psellipse[linestyle=dashed](-50,60)(50,70)

\psline(-55,55)(-20,90)
\psline(-20,80)(-55,115)

\psline(-55,65)(-20,30)
\psline(-20,40)(-55,5)

\uput[0](-45,13){\small $p_V$}

\pscurve(-55,25)(-35,-10)(-55,-45)

\uput[0](-35,-20){\small $\Gamma_V$}

\rput[lt]{0}(20,100){
\parbox{1.9in}{\small $p_a(V) = 0$, $V^2 = -1$ and $V\cdot M = 0$,
where $M$ are the sections given by the marked points $Y_t\cap \pi^*(W)$.}
}

\psline(-180,-10)(-90,-10)
\psline(-180,-20)(-90,-20)
\psline(-180,-30)(-90,-30)
\uput[180](-180,-20){\small $M$}

\uput[180](-105,100){\small $V$}

\end{pspicture}

\caption{Configuration of $V$}
\label{FIG004}
\end{figure}

Let $V\subset Z$ be a maximal connected component of $Z$ and
let us assume that
\begin{equation}\label{E049}
\sum_{\Gamma\subset V} (2p_a(\Gamma) - 2 + |\sigma(\Gamma)|) < 0;
\end{equation}
otherwise, there is nothing to prove. This happens only if (see
\figref{FIG004})
\begin{enumerate}
\item every component of $V$ is smooth and rational;
\item the dual graph of $V$ is a tree;
\item $V$ meets the rest of $Y_0$ at a single point $p_V$ and let us
assume that $p_V = V\cap \Gamma_V$ for some irreducible component
$\Gamma_V \not\subset Z$;
\item $\pi_* V \ne 0$, i.e., $V$ is noncontractible under $\pi$;
\item $V\cap M = \emptyset$, where $M$ is the closure
of $\{Y_t\cap \pi^*(W)\}_{t\ne 0}$ in $Y$;
\item the LHS of \eqref{E049} is exactly $-1$.
\end{enumerate}


Note that $p_V\in \sigma(\Gamma_V)$. If $p_V\not\in \rho(\Gamma_V)$, then
\begin{equation}\label{E050}
|\sigma(\Gamma_V)| - |\rho(\Gamma_V)|
+ \sum_{\Gamma\subset V} (2p_a(\Gamma) - 2 + |\sigma(\Gamma)|) \ge 0.
\end{equation}
If this holds for every maximal connected component $V\subset Z$,
\eqref{E406} easily follows. Otherwise, suppose that
\begin{equation}\label{E051}
p_V\in \rho(\Gamma_V).
\end{equation}
for some $V$. Since $\pi(V)\cap \partial D_J = \emptyset$,
\eqref{E051} can happen only if there exists an irreducible component
$B$ of $W$ such that $\pi(\Gamma_V)\not\subset B$ and
$F_V\subset B$ for some $F_V\subset \pi_* V$. Let us consider the
divisor $\pi^*(B)$ on $Y$. We can write it as
\begin{equation}\label{E052}
\pi^*(B) = \hat M + \hat V + N
\end{equation}
where $\hat M\subset M$, $\supp(\hat V) \subset V$, $F_V\subset \pi(\hat
V)$, $\supp(N)\subset Z_0$ and $\hat V\cap N = \emptyset$. 
Since $V\cap M = \emptyset$, $\hat V\cap \hat M = \emptyset$.
This leads to the contradiction:
\begin{equation}\label{E053}
(\hat V)^2 = \hat V(\hat M + \hat V + N) = \pi_* \hat V\cdot B
\ge 0.
\end{equation}
Therefore, \eqref{E051} cannot hold for any maximal connected
component $V$ of $Z$ and \eqref{E406} follows. This justifies
\eqref{E213} with $\Gamma\subset U_Y$ excluded in the summation.

Now let us turn to \eqref{E215}. Let $G = f(D_\alpha)\subset D_\beta$
for some $\beta\in I$ and $\beta \ne \alpha$. By our assumption on
$f$, $G\cap D_i = \emptyset$ for $i\ne \alpha, \beta$ and $X$ is the
blowup of $\hat X$ along $G$.

Let $\hat Y$ be a flat family of curves that fits in the commutative
diagram:
\begin{equation}\label{E054}
\begin{CD}
Y @>\pi>> X\\
@VV{\phi}V @VVfV\\
\hat Y @>\hat\pi>> \hat X
\end{CD}
\end{equation}
and we assume that $\hat\pi: \hat Y\to \hat X$ is a family of
semistable maps with marked points $\hat Y_t \cap \hat\pi^*(\hat W)$.
Obviously, the map $\phi: Y\to \hat Y$ simply contracts some $-1$ or
$-2$ smooth rational curves $\Gamma$ with $f_* \pi_* \Gamma = 0$.

For an irreducible component $\hat \Gamma\subset \hat Y_0$ with $\hat
\pi_* \hat \Gamma \ne 0$, there exists a unique component
$\Gamma\subset Y_0$ that dominates $\hat\Gamma$ by $\phi$. Actually,
$\phi: \Gamma\to \hat \Gamma$ is a partial normalization of
$\hat\Gamma$. 
Let $\sigma(\Gamma),\rho(\Gamma) \subset \Gamma$ and $\sigma(\hat
\Gamma), \rho(\hat\Gamma)\subset \hat\Gamma$ be the sets defined as
before.

In order to show \eqref{E215}, it is enough to show
\begin{equation}\label{E055}
2 p_a(\hat \Gamma) - 2 + |\sigma(\hat \Gamma)| \ge
2 g(\Gamma) - 2 + |\rho(\Gamma)|
\end{equation}
for every component $\Gamma\subset Y_0$ with $F = \supp(\pi_* \Gamma)\subset
D_\alpha$ and $f_*\pi_*\Gamma\ne 0$, since 
the RHS of \eqref{E055} is at least $(\deg \pi_\Gamma) \Phi_{X,W}(F)$
by \eqref{E405}.

In order to show \eqref{E055}, it in turns suffices to show
\begin{equation}\label{E056}
\rho(\Gamma) \subset \phi^{-1}(\sigma(\hat\Gamma)).
\end{equation}
So it remains to verify \eqref{E056}. Let $p \in \rho(\Gamma)$.

There are two possibilities that $p\in I_{D_\alpha}(\Gamma, G)$ or
$p\in I_X(\Gamma, W)$.

Suppose that $p\in I_{D_\alpha}(\Gamma, G)$. Then by \lemref{LEM005}, there
exists a component $\Gamma'\subset Y_0$ such that
$\pi(\Gamma')\subset D_\beta$, $\pi_* \Gamma' \ne 0$ and $\Gamma$ and
$\Gamma'$ are joined by a chain of curves which are contracted to
$\pi(p)$ by $\pi$ (see \figref{FIG005}). Obviously,
$\phi_* \Gamma' \ne 0$ and hence $p\in \phi^{-1}(\sigma(\hat\Gamma))$.


\begin{figure}[ht]
\centering
\begin{pspicture}(-160,-310)(300,-125)


\uput[45](100,-140){\small $D_\alpha$}
\uput[225](-150,-160){\small $X_0$}

\psline(-100,-140)(100,-140)
\psline(-100,-230)(100,-230)

\psline(200,-200)(180,-130)
\uput[0](190,-165){\small $\Gamma$}
\uput[0](200,-190){\small $p$}

\pscurve(200,-225)(190,-205)(200,-185)
\pscurve(200,-250)(190,-230)(200,-210)
\pscurve(200,-270)(190,-255)(200,-235)

\psline(0,-230)(-20,-160)
\uput[0](-10,-195){\small $F$}

\psline{->}(160,-230)(130,-230)
\uput[90](145,-230){\small $\pi$}

\uput[45](0,-230){\small $\pi(p)$}

\psline(200,-255)(160,-295)
\uput[315](180,-275){\small $\Gamma'$}

\psline(0,-230)(-40,-270)
\uput[315](-20,-250){\small $F'$}

\psline[linestyle=dashed](-50,-180)(100,-180)
\psline(100,-180)(150,-180)
\psline(-150,-280)(50,-280)

\uput[225](-150,-280){\small $D_\beta$}

\uput[270](65,-230){\small $G$}

\psline[linestyle=dashed](-50,-180)(-100,-230)
\psline(-100,-230)(-150,-280)

\psline(-100,-140)(-100,-230)
\psline(100,-140)(100,-230)

\psline(150,-180)(50,-280)

\end{pspicture}

\caption{Application of \lemref{LEM005}}\label{FIG005}
\end{figure}

Suppose that $p\in I_X(\Gamma, W)$ and
$p\not\in \phi^{-1}(\sigma(\hat\Gamma))$. This can happen only if
there exists a connected union of components $V\subset Y_0$
such that
\begin{enumerate}
\item $p = V\cap \Gamma$ and $p$ is the only intersection between $V$
and the rest of $Y_0$;
\item $\phi_* V = 0$ and $V$ is a tree of smooth rational curves;
\item $L = \supp(\pi_* V)$ is a fiber of the projection $f:
D_\alpha\to G$;
\item $V\cap M = \emptyset$;
\item there exists an irreducible component $B\subset W$ such that
$\pi(\Gamma) \not\subset B$ and $L\subset B$.
\end{enumerate}
Again we have the configuration of $V$ as in \figref{FIG004}. Next, we
argue in exactly the same way as for \eqref{E051} by
observing that $B\cdot L \ge 0$. This finishes the proof of 
\eqref{E056} and hence the proof of \eqref{E215}.

By an almost identical argument as before, which we will not repeat,
we can prove
\eqref{E215} with $\Gamma\subset U_Y$ excluded in the summation.

\end{document}